\documentclass[11pt,a4paper]{amsart}
\usepackage[arrow, matrix, curve]{xy}
\usepackage{amssymb}
\usepackage{color}
\usepackage{hyperref}

\input{pstricks.tex}

\usepackage{tikz}
\newtheorem{thm}{Theorem}
\newtheorem{lem}{Lemma}
\newtheorem{cor}{Corollary}

\newtheorem{prop}{Proposition}
\newtheorem{defn}{Definition}

\newenvironment{prf}{{\it Proof:}\quad}{\hfill$QED$}

\begin{document}

\title{Foliated norms on fundamental group and homology}
 \author{Thilo Kuessner}

\maketitle

\begin{abstract}

This paper compares two invariants of foliated manifolds which
seem to measure the non-Hausdorffness of the leaf space: the transversal
length on the fundamental group and the foliated Gromov norm on the homology.
We consider foliations with the property that the set of singular 
simplices strongly
transverse to the foliation satisfies a weakened version of the Kan
extension
property. (We prove that this assumption is fairly general: it holds for
all fibration-covered foliations, in particular
for all foliations of 3-manifolds without Reeb components.)
For such foliations we show that
vanishing of the transversal length implies triviality of the foliated
Gromov norm, and, more generally, that uniform bounds on the transversal
length imply explicit bounds for the foliated Gromov norm. This is
somewhat surprising in view of the fact that transversal length is defined
in terms of 1- and 2-dimensional objects.\end{abstract}

%{\bf Introduction}\\

Foliations are one of the tools to study the topology of 3-manifolds. Let $M$ be a closed 3-manifold and
$\mathcal{F}$ a $C^\infty$-foliation without Reeb components. It is known that the pullback foliation $\widetilde{\mathcal{F}}$ of the universal covering $\widetilde{M}$ is a foliation of ${\Bbb R}^3$ by leaves homeomorphic to
${\Bbb R}^2$. By Palmeira's theorem, such foliations $\widetilde{\mathcal{F}}$
are completely classified by their leaf space, a simply connected non-Hausdorff
1-manifold. Therefore, to study codimension one foliations without Reeb components on 
3-manifolds, it is useful
to study the actions of 3-manifold fundamental groups on simply connected 1-manifolds. In 
particular, a foliation can be complicated in two ways: 
the leaf space of the pull-back foliation on $\widetilde{M}$ 
can be very branched (i.e., be non-Hausdorff with a complicated 
pattern of branching points), or the group acting on this 1-manifold 
can be complicated. The second point is of course inherent in $M$ itself: the acting group is $\pi_1M$, which does not depend
on the specific foliation ${\mathcal F}$. Thus, an invariant measuring the complexity of $\mathcal{F}$ should be composed of $\pi_1M$ and an invariant describing the branching of the leaf space of $\widetilde{\mathcal{F}}$.

Calegari defined in \cite{cale} an invariant $\parallel M\parallel_{\mathcal{F}}$, the foliated Gromov norm, and proved several theorems showing that the size of this invariant is related to the branching of the leaf space. 
His invariant is a refinement of the Gromov norm $\parallel M\parallel$, which measures the complexity of a manifold $M$ and has subtle relations with the fundamental group $\pi_1M$. The foliated Gromov norm $\parallel M\parallel_{\mathcal{F}}$ is at least as large as $\parallel M\parallel$ and it is the difference $\parallel M
\parallel_{\mathcal{F}} -
\parallel M\parallel$ which seems to be related to the branching of the leaf space of
$\widetilde{\mathcal{F}}$. In particular, Calegari proved that $\parallel 
M
\parallel_{\mathcal{F}}=\parallel 
M
\parallel$ if the leaf space does not
branch or branches in only one direction, and he exhibited large classes of branching foliations where the foliated Gromov norm is strictly larger than
the simplicial volume.

We consider a second invariant $l_{\mathcal{F}}$, which is a pseudonorm
on the fundamental group of a foliated manifold. (The definition is reminiscent of a similar definition in \cite{cale}.) We show that, under a technical assumption, vanishing of $\l_{\mathcal{F}}$ implies triviality of the foliated Gromov norm and, more generally,
the foliated Gromov norm can be bounded in terms of $l_{\mathcal{F}}$ and
$\parallel M\parallel$.\\
\\
To describe the technical assumption, we need to sketch two definitions (which will be made precise in the first chapter). A singular simplex is said to be strongly transverse to $\mathcal{F}$ if the induced foliation is affine and there is no 'backtracking' (see Section 1.1). We say that $\mathcal{F}$ satisfies the weak Kan property in degrees $\ge 2$ if, for 
$2\le n\le \dim\left(M\right)$, any (n+1)-tuple of strongly
transversal n-simplices with compatible boundaries
admits an (n+1)-simplex, with the given simplices as boundary faces,
whose (n+2)-th boundary face is strongly transverse as well (see Section 1.2). \\
Moreover we will need the condition that $\widetilde{M}/\widetilde{\mathcal{F}}$ is a 1-connected 1-manifold. This was stated as a theorem in chapter 1 of \cite{palm}, but an inspection of the proof shows that it assumes that 
$\widetilde{\mathcal{F}}$ is a foliation of ${\Bbb R}^n$ by hyperplanes. In particular, the condition is satisfied
for Reebless foliations of 3-manifolds.\\
For foliations satisfying these conditions we have:\\
\\
{\bf Corollary 4:} {\em 
Let $\left(M,{\mathcal{F}}\right)$ be a foliated manifold such that
$\widetilde{M}/\widetilde{\mathcal{F}}$ is a 1-connected 1-manifold
and the set of 
strongly transversal simplices satisfies the weak Kan property
in degrees $\ge 2$. Then 
$$\left(\exists x_0\in M\forall \gamma\in\pi_1\left(M,x_0\right): l_{\mathcal{F}}\left(\gamma\right)=0\right)\Longrightarrow\parallel M\parallel_{\mathcal{F}}=
\parallel M\parallel.$$}\\
\\
More generally, we show 
that, under the same assumptions, 
the nontriviality of the
foliated Gromov norm can be estimated in terms of the norm $l_{\mathcal{F}}$,
if the latter happens to be uniformly bounded.\\
\\
{\bf Theorem 1:} {\em Let $\left(M,{\mathcal{F}}\right)$ be a foliated manifold such that
$\widetilde{M}/\widetilde{\mathcal{F}}$ is a 1-connected 1-manifold
and the set of strongly transversal simplices satisfies the weak Kan property
in degrees $\ge 2$. Then}
$$\parallel M\parallel_{\mathcal{F}}\le
\left(1+\left(\dim\left(M\right)+1\right)sup\left\{l_{\mathcal{F}}\left(\gamma\right)
:\gamma\in\pi_1
\left(M,x_0\right)\right\}\right)\parallel M\parallel.$$
\\
The question of interest is then, of course,
for which foliations the weak Kan property in degrees $\ge 2$ holds. Let
$K^{\mathcal{F}}_{st}$ denote the set of strongly transversal simplices 
and $K^{\mathcal{F}}$
the set of transversal simplices.
We show that a fairly general class of foliations satisfies the weak Kan property in degrees $\ge2$, in particular\\
\\
{\bf Lemma 4:} {\em Let $M^3$ be a compact 3-manifold, $\mathcal{F}$ a codimension one foliation
on $M$. Then \\
a) the simplicial sets $K_{st}^{\mathcal{F}}$ and $K^{\mathcal{F}}$
satisfy the weak Kan property in degree 3,\\
b) the simplicial set $K^{\mathcal{F}}$
does not satisfy the weak Kan property in degree 2, \\
c) the simplicial set $K_{st}^{\mathcal{F}}$ satisfies
the weak Kan property in degree 2 if and only if $\mathcal{F}$ has no Reeb component
.}\\
\\
More generally, we have\\
\\
{\bf Theorem 2:} {\em 
If $M$ is an $m$-dimensional
manifold and $\mathcal{F}$ a fibration-covered codimension
one foliation, such that $\widetilde{M}/\widetilde{\mathcal{F}}$ is a 1-connected 1-manifold, then 

a) $K_{st}^{\mathcal{F}}$ satisfies
the weak Kan property in degrees $2,\ldots,m$, and

b) $K^{\mathcal{F}}$ satisfies the weak Kan property in degrees $3,\ldots,m$.}\\
\\
{\bf Corollary 6:} {\em 
If $M$ is an $m$-dimensional 
manifold and $\mathcal{F}$ a codimension 
one foliation satisfying the following conditions:\\
(i) for every leaf $F$ is $\pi_1F\rightarrow\pi_1M$ injective,\\
(ii) for every leaf $F$, the universal covering 
$\widetilde{F}$ is homeomorphic to ${\Bbb R}^{m-1}$,\\
then

a) $K_{st}^{\mathcal{F}}$ satisfies
the weak Kan property in degrees $2,\ldots,m$, and

b) $K^{\mathcal{F}}$ satisfies the weak Kan property in degrees $3,\ldots,m$.}\\
\\
It seems worth mentioning that the assumption in theorem 10 is not just technical but can actually not be avoided. For example, it follows from \cite[Theorem 2.6.2]{cale} 
that there are foliations ${\mathcal{F}}$
of the 3-sphere ${\Bbb S}^3$ with arbitrarily large Gromov norm $\parallel {\Bbb S}^3\parallel_{\mathcal{F}}$, but $\parallel{\Bbb S}^3\parallel=0$ and $l_{\mathcal{F}}=0$ since $\pi_1{\Bbb S}^3=0$. Thus an inequality as in Theorem 1 can not hold true for foliations with Reeb
components.\\

Our results and proofs suggest that the branching of the leaf space is directly related to 
the failure of the Kan extension property. (For example, the nontriviality of $l_{\mathcal{F}}$ expresses the failure of the weak Kan property in degree 1.) Thus it would be nice to have quantitative invariants of
sets of simplices which measure the failure of the Kan extension property, as this would give interesting invariants of foliations.

Conventions: We assume all manifolds and foliations to be $C^\infty$,
tangentially and transversally. (This will be needed for applications of Palmeira's
theorems.)
We assume all manifolds to be orientable. All theorems generalize in an obvious way to non-orientable manifolds.\\
I thank the referee for some helpful comments.

\section{Preparations}

\subsection{Basic definitions}

Let $M$ be a manifold and $\mathcal{F}$ a {\bf codimension one foliation} of $M$.
Let $\Delta^n$ be the standard simplex in ${\Bbb R}^{n+1}$, and $\sigma:\Delta^n\rightarrow M$ some (continuous, not necessarily differentiable) singular simplex.

The foliation $\mathcal{F}$ induces an equivalence relation on $\Delta^n$ by:
$x\sim y\Longleftrightarrow \sigma\left(x\right)$ and $\sigma\left(y\right)$ belong to the same connected component of $L\cap\sigma\left(\Delta^n\right)$ for some leaf $L$ of $\mathcal{F}$.\\
This equivalence relation may or may not be induced by some foliation of the standard simplex $\Delta^n$. 

We say that a singular simplex $\sigma:\Delta^n\rightarrow M$ is {\bf foliated} if the equivalence relation $\sim$ is induced by a foliation of $\Delta^n$. We will denote this foliation of $\Delta^n$ by ${\mathcal{F}}\mid_{\sigma}$.\\

We call a foliation $\mathcal{F}$ of $\Delta^n$ affine if there is an affine mapping $f:\Delta^n\rightarrow{\Bbb R}$ such that $x,y\in\Delta^n$ belong to the same leaf if and only if $f\left(x\right)=f\left(y\right)$. We say that a foliation of $\Delta^n$ is conjugate to an affine foliation if there is a simplicial homeomorphism $H:\Delta^n\rightarrow\Delta^n$ such that $H^*\mathcal{F}$ is an affine foliation.

We say that a singular $n$-simplex $\sigma:\Delta^n\rightarrow M$, $n\ge 2$, is {\bf transverse} to $\mathcal{F}$ if it is foliated and  it is \\
- either contained in a leaf,\\
- or the induced foliation ${\mathcal{F}}\mid_\sigma$ is conjugate to an affine foliation $\mathcal{G}$ of $\Delta^n$.\\
\\
For $n=1$, we say that a singular 1-simplex $\sigma:\Delta^1\rightarrow M$ is transverse to $\mathcal{F}$ if it is\\
- either contained in a leaf,\\
- or for each foliation chart $\phi:U\rightarrow {\Bbb R}^{m-1}\times{\Bbb R}^1$ (with m-th coordinate map 
$\phi_m:U\rightarrow{\Bbb R}^1$) one has that $\phi_m\circ\sigma\mid_{\sigma^{-1}\left(U\right)}:
\sigma^{-1}\left(U\right)\rightarrow{\Bbb R}^1$ is locally surjective at 
all points of $int\left(\Delta^1\right)$, i.e.\ for all 
$p\in int\left(\Delta^1\right)\cap \sigma^{-1}\left(U\right)$, the 
image of $\phi_m\circ\sigma\mid_{\sigma^{-1}\left(U\right)}$ contains a neighborhood of $\phi_m\circ\sigma\left(p\right)$.\\
\\
We say that a singular simplex $\sigma:\Delta^n\rightarrow M$ is
{\bf strongly transverse} if it is transverse and \\
- either contained in a leaf,\\
- or for each foliation chart $\phi:U\rightarrow {\Bbb R}^{m-1}\times{\Bbb R}^1$ (with m-th coordinate map
$\phi_m:U\rightarrow{\Bbb R}^1$) one has that $\phi_m\circ\sigma\mid_{\sigma^{-1}\left(U\right)}:
\sigma^{-1}\left(U\right)\rightarrow{\Bbb R}^1$ is locally surjective at all
points of $int\left(\Delta^n\right)$.\\
(In particular, for $n=1$ transversality and strong transversality are the same.)

In general, if we have a preferred set of 'transversal' simplices $T$, we define the transversal Gromov norm of a compact, orientable manifold $M$ with fundamental class $\left[M,\partial M\right]\in H_n\left(M,\partial M;{\Bbb R}\right)$ as $$\parallel M,
\partial M\parallel_T:=inf\left\{\sum \mid a_i\mid: \sum a_i\sigma_i \mbox{ represents }\left[M,\partial M\right],
\sigma_i\in T\right\}.$$
In the case of a codimension one foliation $\mathcal{F}$ on $M$, taking for $T$ the set $K^{\mathcal{F}}$
of singular simplices transverse to $\mathcal{F}$ (as defined above) and defining $$\parallel M,\partial M\parallel_{\mathcal{F}}:=\parallel M,\partial M\parallel_{K^{\mathcal{F}}},$$ one gets Calegari's definition of foliated Gromov norm (\cite{cale}).\\
If $\partial M=\emptyset$, we will omit $\partial M$ from the notation.\\
We will mainly work not with $\parallel M,\partial M\parallel_{\mathcal{F}}$, but with
the following notion: let $K_{st}^{\mathcal{F}}$ be the set of simplices strongly transverse to $\mathcal{F}$
and define $$\parallel M,\partial M\parallel_{\mathcal{F}}^{st}:=\parallel
M,\partial M\parallel_{K_{st}^{\mathcal{F}}}.$$

There is an obvious inequality 
$\parallel M,\partial M\parallel_{\mathcal{F}}\le\parallel M,\partial M\parallel_{\mathcal{F}}^{st}$.

\subsection{Weak Kan property and Gromov norm}

Let $X$ be a topological space and $S_k\left(X\right)$ its set of singular 
$k$-simplices
with the face maps $\partial_i:S_k\left(X\right)\rightarrow
S_{k-1}\left(X\right)$ and the degeneration maps $s_i:S_k\left(X\right)\rightarrow
S_{k+1}\left(X\right)$ for $0\le i\le k$. 

We define in this section some properties of subsets $K\subset S_*\left(X\right)$. The cases we are interested in 
is that $X=M$ is a foliated manifold, with fixed foliation $\mathcal{F}$, and $K$ is either $K^{\mathcal{F}}$ or $K_{st}^{\mathcal{F}}$. 
It is easy to see that $K_{st}^{\mathcal{F}}$ is closed under face 
and degeneration maps. On the other hand,
$K^{\mathcal{F}}$ is not a simplicial set. 
Indeed, consider a 1-simplex which is not transverse. Then its degeneration 
$s_0\left(\sigma\right)$ is transverse but $\partial_0s_0\left(\sigma\right)=\sigma$ 
is not. Thus $\partial$ does not map $K_2^{\mathcal{F}}$ to $K_1^{\mathcal{F}}$. (However, $K^{\mathcal{F}}$ is clearly closed under
degeneration maps.)\\
Although we are concerned with $K^{\mathcal{F}}_{st}$, we will for completeness also discuss the properties of $K^{\mathcal{F}}$ in this paper. Therefore we do not assume that $K$ is a simplicial subset.

A set $T\subset S_*\left(X\right)$ is called a Kan complex if it is a simplicial
set, i.e.\ stable with respect to face and degeneration maps, and 
if the following holds:\\
for any collection of $n+1$ $n$-simplices
$\left\{
\tau_0,\ldots,\tau_{k-1},\tau_{k+1},\ldots,\tau_{n+1}\right\}\subset T_n$ with
$\partial_i\tau_j=\partial_{j-1}\tau_i$ for all $i<j$,
there exists an $n+1$-simplex $\sigma\in T_{n+1}$ with $\partial_i\sigma=\tau_i$ 
for all $i\not=k$.\\

The theory of Kan complexes is well developed. However, we will need to work with 
sets of simplices which only satisfy the following condition.\\

{\bf Weak Kan property.} We say that a set $T\subset 
S_*\left(X\right)=\cup_{k\in {\Bbb N}}S_k\left(X\right)$, satisfies the weak Kan property in degree n if the following conditions hold: \\
- for any collection of $n+1$ $n$-simplices
$\left\{
\tau_0,\ldots,\tau_{k-1},\tau_{k+1},\ldots,\tau_{n+1}\right\}\subset T_n:=T\cap S_n\left(X\right)$ with 
$\partial_i\tau_j=\partial_{j-1}\tau_i$ for all $i<j$,
there exists an $n+1$-simplex $\sigma\in S_{n+1}\left(X\right)$
with $\partial_i\sigma=\tau_i$ for
all $i\not=k$, such that $\partial_k\sigma\in T_n$ (not necessarily $\sigma\in T_{n+1}$). \\

Moreover, we will need the following notion: let $\Delta^n,\Delta^{n-1}$ be standard simplices, $r:\Delta^n\rightarrow\Delta^{n-1}$ be any affine 
mapping with $r\left(v_i\right)=v_i$ for 
$0\le i\le n-1$ and $r\left(v_n\right)\in\Delta^{n-1}$ arbitrary, then 
we say that a
singular n-simplex $\sigma:\Delta^n\rightarrow X$ is a {\bf general degeneration} of a singular simplex $\tau:\Delta^{n-1}\rightarrow X$ if $\sigma=r\tau$. \\
We say that $K\subset S_*\left(X\right)$ is closed under general degenerations
if $\tau\in K_{n-1}$ implies $r\tau\in K_n$ for any such $r$.\\
\\
{\bf Simplicial approximation property.} For a singular simplex $\sigma:\Delta^n\rightarrow X$ and $k\in{\Bbb N}$ let $sd^k\left(\sigma\right):sd^k\left(\Delta^n\right)\rightarrow X$ denote
the $k$-th barycentric subdivision 
of $\sigma$, i.e., $sd^k\left(\sigma\right)=\sigma\Phi^k$ where $\Phi^k:sd^k\left(\Delta^n\right)\rightarrow\Delta^n$
is the canonical continuous projection.
We say that a subset $K\subset S_*\left(X\right)$ satisfies the simplicial approximation property if the following holds true:\\
for each $j\in {\Bbb N}$ and each singular simplex $\sigma\in S_*\left(X\right)$ with $j$-skeleton in $K$, there is $k\in{\Bbb N}$ such that $sd^k\left(\sigma\right)$ is homotopic to a simplicial mapping $f:sd^k\left(\Delta^n\right)\rightarrow K$, by a homotopy which {\em leaves the $j$-skeleton of $\sigma$ pointwise fixed}.\\
\\
{\bf Geometric Realisation.} For a subset $K\subset S_*\left(X\right)$, which is closed under degeneration maps, define 
its geometric realisation $RK$
exactly as in \cite{curt}, p.118, except for the following: \\
if, for some simplex $x\in K$, some boundary face $\partial_ix$ does
not belong to $K$, then we erase this (open) boundary face from the image of $x$
in $RK$ (but we do not erase iterated boundaries of $\partial_ix$ in case they belong to $K$). That is, the image of $x$ in $RK$ will not necessarily be closed. 
(In other words, we consider $RK$ as a subset $RS_*\left(X\right)$, with $RS_*\left(X\right)$ defined in \cite{curt}, such that a point in $RS_*\left(X\right)$ belongs to $RK$ if it is in the image of some simplex $x\in K$.)\\
\\
{\bf Simplicial and singular Gromov norm.}
For any subset $K\subset S_*\left(X\right)$, which is closed under degeneration maps, 
there is its geometric realisation $RK$. Even though $K$ need not be a simplicial set, 
we can define define its group of cycles $Z_*\left(K;{\Bbb R}\right):=
\left\{\sum_{i=1}^r a_i\sigma_i: a_i\in{\Bbb R},\sigma_i\in K, \partial\sum_{i=1}^r 
a_i\sigma_i=0\right\}$ and its group of boundaries $B_*\left(K;{\Bbb R}\right)$ as boundaries 
of chains in $K$, and thus its homology $H_*\left(K;{\Bbb R}\right)=Z_*
\left(K;{\Bbb R}\right)/B_*\left(K;{\Bbb R}\right)$. The inclusion $i:K\rightarrow S_*^{sing}\left(RK\right)$ induces then a homomorphism
$$i_*:H_*\left(K;{\Bbb R}\right)\rightarrow H_*^{sing}\left(RK;{\Bbb R}\right)$$ 
of the simplicial homology of $K$ to the singular homology of $RK$.\\
Note that $i_*$ factors over $\Phi_{K*}$, where $\Phi_K:K\rightarrow S_*^{simp}\left(RK\right)$ is the canonical mapping.\\
For a homology class $h\in H_n\left(K;{\Bbb R}\right)$ define $$\parallel h\parallel_{simp}:=inf\left\{\sum \mid a_i\mid: \sum a_i\sigma_i \mbox{ represents }h\right\}$$ 
and $$\parallel h\parallel_{sing}:=inf\left\{\sum\mid a_i\mid:\sum a_i\sigma_i\mbox{ represents }\Phi_{K*}h\in H_*^{sing}\left(RK;{\Bbb R}\right), \sigma_i 
\mbox{ singular simplices}\right\}.$$ Of course, $\parallel h\parallel_{sing}\le\parallel h\parallel_{simp}$.

We hope that it does not lead to confusion that the well-known Gromov norm on a
topological space $X$ is the {\bf simplicial} Gromov norm of the simplicial set $K=S_*\left(X\right)$, meanwhile the singular 
Gromov norm on $S_*\left(X\right)$ might be smaller.\\
\\

To motivate the results of this subsection, we first discuss the case of {\bf Kan complexes}.

For a simplicial set $K$, let $H^*_{b,simp}\left(K\right)$ be the simplicial 
bounded cohomology of $K$ and $H^*_{b,sing}\left(RK\right)$ be the
singular bounded cohomology of the geometric realization (see \cite{grom} for 
the definition of bounded cohomology). In general, of course, $H^*_{b,simp}$
and $H^*_{b,sing}$ are unrelated. (For example, if $K$ is a finite 
simplicial complex, then 
$H^*_{b,simp}\left(K\right)=H^*\left(K\right)$ but, in general,
$H^*_{b,sing}\left(RK\right)\not=H^*\left(RK\right)$.) However, for Kan complexes, there is the following lemma (which is similar in spirit, but not directly related, to the isometry lemma, \cite{grom}, p.43).
\begin{lem}\label{kanc} If $K$ is a Kan complex, then \\
a) $H^*_{b,simp}\left(K\right)$ is
isometrically isomorphic to $H^*_{b,sing}\left(RK\right)$,\\
b) $\parallel h\parallel_{sing}=\parallel h\parallel_{simp}$ for all $h\in H_*\left(K\right)$.\end{lem}
\begin{prf}
Let $\Psi_{RK}:RS_*\left(RK\right)
\rightarrow RK$ be the canonical continuous mapping which projects each singular simplex to its image.
By the simplicial extension theorem (which is proved in \cite{curt}, where it is attributed to unpublished work of Barratt and Kan), there exists a 
simplicial mapping $g:S_*\left(RK\right)\rightarrow K$ such that $R\left(g\right)$ is
homotopic to $\Psi_{RK}$.

Let $\Phi_K:K\rightarrow S_*\left(RK\right)$ be the 
canonical simplicial mapping. It is shown in \cite{curt} that $\Phi_Kg$ and $g\Phi_K$ are chain homotopic to the identities. By dualizing we get isometric isomorphisms of bounded cohomology. ($g_*$ and $\Phi_{K*}$ 
do not increase norms, thus they must be isometries since their composition is the identity.) This proves claim a.
By the well-known duality between the norm in bounded cohomology and the Gromov norm on homology (see \cite{grom}), claim b follows.\end{prf}\\
\\
We come now to a general result about weak Kan complexes which will be needed for the proof of \hyperref[ineq]{Theorem \ref*{ineq}}.
\begin{prop}\label{sap} Let $X$ be a topological space, $n\in{\Bbb N}$, and
$K\subset S_*\left(X\right)$
a set of singular simplices which 
satisfies the simplicial approximation property,
satisfies the weak Kan property in degrees $\ge n$ and is closed with respect to 
general degenerations.
Let $h\in H_*\left(K;{\Bbb R}\right)$ be a homology class and let $\sum_{j=0}^sa_j\sigma_j\in C_*\left(X;{\Bbb R}\right)$ be a cycle such that\\
i) the n-1-skeleta of $\sigma_1,\ldots,\sigma_s$ belong to $K$\\
ii) $\sum_{j=0}^s a_j\sigma_j$ represents $i_*h\in H_*\left(X;{\Bbb R}\right)$.\\
Then $\parallel h\parallel_{simp}\le\sum_{j=0}^s\mid a_j\mid$.\end{prop}
A special case (with $X=RK$) is the following corollary.
\begin{cor}\label{sap1}
Let $Y$ be a topological space, $n\in{\Bbb N}$, and $K\subset S_*\left(Y\right)$ a set of singular simplices which satisfies 
the simplicial approximation property, 
satisfies the weak Kan property in degrees 
$\ge n$,
is closed with respect to general degenerations,
and which is such that each cycle $\sum_{j=0}^s a_j\sigma_j\in C_*^{sing}\left(RK;{\Bbb R}\right)$ is homologous to a cycle $\sum_{j=0}^sa_j\sigma_j^\prime\in C_*^{sing}\left(RK;{\Bbb R}\right)$ (with equal coefficients)
with the n-1-skeleta of $\sigma_1^\prime,\ldots,\sigma_s^\prime$ belonging to $K$. Then $\parallel h\parallel_{sing}=\parallel h\parallel_{simp}$ for all $h\in H_*\left(K\right)$.\end{cor}

{\em Proof of \hyperref[sap]{Proposition \ref*{sap}}:}\\
Let $L$ be the simplicial set built of $\sigma_1,\ldots,\sigma_s$ together with
all of their (iterated) faces
and degenerations. Let $p:RL\rightarrow RS_*\left(X\right)$ be the canonical continuous mapping
which projects each $\sigma$ to its image. By assumption (i), $p$ maps the n-1-skeleton
of $L$ to $RK\subset RS_*\left(X\right)$.

By the simplicial approximation property, there exists some $k\in{\Bbb N}$ and some simplicial mapping $f^\prime:sd^k\left(L\right)\rightarrow K$ such that $R\left(f^\prime\right)\sim p\Phi^k$ and $R\left(f^\prime\right)=p\Phi^k$ on
$\left(\Phi^k\right)^{-1}\left(L_{n-1}\right)$. (We may choose $k$ uniformly
because $L$ contains only finitely many nondegenerate simplices.)

As in \cite{curt}, chapter 12, 
we are looking for a simplicial mapping $g:L\rightarrow K$ and a simplicial mapping $F:sd^k\left(L\right)\times I\rightarrow K$ such that $R\left(F\right)$
provides a homotopy between $R\left(f^\prime\right)$ and $R\left(g\right)\Phi^k$.
Since $S_*\left(X\right)$ is a Kan complex, one can apply the construction
in \cite{curt} to construct such mappings with image in $S_*\left(X\right)$,
not necessarily in $K$.
Our task is to prove that (under the assumptions of proposition 1), we can construct $F$ and $g$ with image in $K$. It suffices to consider the case $k=1$, as the general case follows.

To this aim, we examine the construction in \cite{curt}. There one defined an 
ordering on the simplices of $sd^k\left(x\right)\times I$, 
for each simplex $x\in L$, which allows us
to inductively define $F$ (and finally $g\left(x\right)$), once $F$ was
defined on iterated boundaries of $x$. 

For each simplex $x\in L$, let $W_x$ be the simplices in the canonical triangulation of $sd\left(x\right)\times I$, as used in \cite{cur}, p.204/205. There is a canonical continuous projection 
$\pi:sd\left(x\right)\times
I\rightarrow sd\left(x\right)$ which is not simplicial,
but which actually fails to be simplicial only for one simplex, namely the 'last' simplex
$\overline{w}$ for the ordering of $W_x$ given in \cite{curt}, and for its last boundary face $\partial_{deg\left(x\right)+1}\overline{w}$. (Recall that $g\left(x\right)$ 
was defined to be $g\left(x\right):=\partial_{deg\left(x\right)+1}\overline{w}$.)

We wish to show that the construction in \cite{curt} can be carried out such that, for all $x\in L$ and $w\in W_x$, we have that all boundary faces are mapped to $K$, i.e. that
$F\left( \partial_iw\right)\in K$. This implies especially
$g\left(x\right)=F\left(\partial_{deg\left(x\right)+1}\overline{w}\right)\in K$.

We proceed by induction on the dimension of $x$. Assume we have proved $F\left(\partial_iw_y\right)\in K$ for $\dim\left(y\right)\le m-1$. For the proof of the inductive step we will need to distinguish the cases
$m\le n-1, m=n$ and $m\ge n+1$.

Assume $\dim\left(x\right)=m\le n-1$. We do already know by assumption (i)
that $x$ is mapped to $K$. For each $w\in W_x$, $F\left(w\right)$ 
is a general degeneration of $p\left(x\right)$. Hence $p\left(x\right)\in K$ implies $F\left(w\right)\in K$ for all $w\in W_x$. (In particular,
$g\left(x\right)\in K$ which of course might have been achieved by setting $g\left(x\right)=p\left(x\right)$. We will however need $F\left(w\right)\in K$ for the inductive argument.)

Assume $\dim\left(x\right)=m=n$. Let $x\in W_x$. If either $w\not=\overline{w}$, or $w=\overline{w}$ but $i\not=n+1$, then $F\left(\partial_iw\right)$ is the general degeneration of an n-1-simplex in $K$,
thus $F\left(\partial_iw\right)\in K$. In particular, $F\left(\partial_i\overline{w}\right)\in K$
for $i\not=n+1$. Since $K$ satisfies the weak Kan property in degree $n$, we may choose 
$F\left(\overline{w}\right)$ such that $\partial_iF\left(\overline{w}\right)=
F\left(\partial_i\overline{w}\right)$ for $i\le n$ and $g\left(x\right):=\partial_{n+1}F\left(\overline{w}\right)\in K$.

Assume $\dim\left(x\right)=m\ge n+1$. Then apply the weak Kan property in degree $m$, as in \cite{curt}, to get $F\left(w\right)\in K$ for all $w\in W_x$. This finishes the inductive argument.

Thus we arrive at a simplicial mapping $g:L\rightarrow K$ such that $R\left(F\right)$ provides a homotopy between $R\left(f^\prime\right)$ and $R\left(g\right)\Phi$.
In particular, $g_*\left(\sum_{j=1}^s a_j\sigma_j\right)\in K\subset 
C_*\left(X\right)$ is homologous to $p_*\left(\sum_{j=1}^s a_j\sigma_j\right)$ 
in $C_*\left(X\right)$. This shows \\
$\parallel h\parallel_{simp}\le \sum_{j=1}^s\mid a_j\mid$.
\hfill$\Box$\\
\\
{\bf Relative version.} If $L\subset K\subset S_*\left(X\right)$,
then the geometric realisations satisfy $RL\subset RK$ and
we have a canonical homomorphism $i_*:
H_*\left(RK,RL\right)\rightarrow H_*\left(K,L\right)$ of the simplicial into the singular homology. We define, for $h\in H_*\left(K,L\right)$, the simplicial Gromov norm \\
$\parallel h\parallel_{simp}=inf\left\{\sum\mid a_i\mid: \sum a_i\sigma_i\mbox{ represents }h\right\}$ and the singular Gromov norm \\
$\parallel h\parallel_{sing}=inf\left\{\sum\mid a_i\mid: \sum a_i\sigma_i \mbox{ represents }i_*h\right\}$. A straightforward generalization of the proof of proposition 2 shows: 
\begin{lem}\label{sap2}
Let $X^\prime\subset X$ be topological spaces, $n\in{\Bbb N}$, and
$K\subset S_*\left(X\right)$
a set of singular simplices which
satisfies the simplicial approximation property,
satisfies the weak Kan property in degrees $\ge n$ and is closed with respect to
general degenerations. Let $K^\prime:=K\cap S_*\left(X^\prime\right)$.
Let $h\in H_*\left(K,K^\prime;{\Bbb R}\right)$ be a homology class and let $\sum_{j=0}^sa_j\sigma_j\in C_*\left(X,X^\prime;{\Bbb R}\right)
$ be a relative cycle such that\\
i) the n-1-skeleta of $\sigma_1,\ldots,\sigma_s$ belong to $K$\\
ii) $\sum_{j=0}^s a_j\sigma_j$ represents $i_*h\in H_*^{sing}\left(X,X^\prime;{\Bbb R}\right)$.\\
Then $\parallel h\parallel_{simp}\le\sum_{j=0}^s\mid a_j\mid$.\end{lem}
\noindent
{\bf Relation with \cite{grom}}. Even 
though this is not related to the rest of our paper, we want to 
mention that the framework of Kan complexes can be used to give an alternative proof of results in \cite{grom}, such as the theorem that the fundamental group determines the bounded cohomology. (This was brought to my attention by Elmar Vogt.) We outline the argument.
Let $T=S_*\left(X\right)$. By the construction in \cite{mayp}, p.36, there is a minimal Kan subcomplex $M\subset T$.
Let $\Gamma$ be the set of simplicial automorphisms of $M$,
$\Gamma_n
\subset\Gamma$ the subgroup which fixes the n-skeleton pointwise, and $M_n=M/\Gamma_n$.\\
$M=M_\infty\rightarrow\ldots\rightarrow M_n 
\rightarrow M_{n-1}\rightarrow\ldots M_2\rightarrow M_1$ is the {\bf Postnikov
 system} considered e.g.\ in \cite{mayp},
par.8. Note that $M_1=K\left(\pi_1M,1\right)=K\left(\pi_1T,1\right)$ follows
from \cite{mayp}, thm.8.4. \\
{\em Claim:} $p_n:M_n\rightarrow M_{n-1}$ induces an isometric
isomorphism in bounded cohomology.\\
{\em Outline of proof:}
If $\sigma$ is an n-simplex and $\gamma\in\Gamma_{n-1}$, then
$\sigma$ and $\gamma\sigma$ have the same boundary,
hence define an n-simplex $\sigma *\gamma\sigma$
with boundary in some vertex $v$ of $\sigma$. If $\gamma\in\Gamma_n$, then
$\sigma *\gamma\sigma$ is homotopic to $s_0^n\left(v\right)$.
This defines a map
$I: \Gamma_{n-1}/\Gamma_n\rightarrow\bigoplus_{\sigma\in M_n}\pi_n\left(M
_n,v_\sigma\right)$
which, analogously to \cite{grom},
is an injective group homomorphism. With \cite{mayp}, prop.4.4 this implies that $\Gamma_{n-1}/\Gamma_n$ is abelian,
thus amenable. As in \cite{grom}, one defines $A_n: C_n^b\left(M_n\right)
\rightarrow C_n^b\left(M_{n-1}\right)$ by averaging bounded cochains over
 the orbits of the amenable
group $\Gamma_{n-1}/\Gamma_n$. Clearly, $A_n^*p_n^*=id$ and $\parallel A_n
\parallel\le 1, \parallel p_n\parallel\le 1$, thus $\parallel A_n\parallel
=\parallel p_n\parallel=1$.
It remains\footnote{Note that $\gamma\in\Gamma$ is not homotopic to the 
identity, so the argument in \cite{gro} does not work.} to prove that $p_n
^*A_n^*=id$. \\
Since $p_n:M_n\rightarrow M_{n-1}$ induces an isomorphism of $\pi_1$, we
can lift $p_n$ to a $\pi_1M$-equivariant mapping $\tilde{p}_n:\tilde{M}_n
\rightarrow\tilde{M}_{n-1}$. Let $C_b^n\left(\tilde{M}_n\right)^{\pi_1M}
\simeq C_b^n\left(M_n\right)$ be the $\pi_1M$-invariant
bounded cochains on $\tilde{M}_n$
and let $\tilde{A_n}:C_b^n\left(\tilde{M}_n\right)\rightarrow C_b^n\left(
\tilde{M}_{n-1}\right)$ be defined by averaging over $\Gamma_{n-1}/\Gamma
_n$. The restriction of $\tilde{A}_n$ to the $\pi_1M$-invariant bounded 
cochains
gives $A_n$. Now, $C_b^i\left(\tilde{M}_n\right)$ is a relatively injective
$\pi_1M$-module, and the resolution ${\Bbb R}\rightarrow C_b^0\left(
\tilde{M}_n\right)\rightarrow C_b^1\left(\tilde{M}_n\right)\rightarrow\ldots
$ admits a contracting homotopy, hence, by a standard argument in homological
algebra (\cite{ivan}, p.1099), any two chain maps extending $id_{\Bbb R
}$ are chain homotopic, in particular $\tilde{p}_n\tilde{A}_n\sim id$. 
Restricting to $\pi_1M$-invariant bounded cochains we get $p_nA_n\sim
id$.
This shows that $M_n\rightarrow M_{n-1}$ induces an isometric isomorphism in bounded cohomology and we conclude:
\begin{cor}\label{post} If $S$ and $T$ are Kan complexes and $i:S\rightarrow T$ a simplicial mapping
such that
$i_*:\pi_1S\rightarrow\pi_1T$ is an isomorphism, then $i^*:
H_b^*\left(T\right)\rightarrow H_b^*\left(S\right)$ is an
isometric isomorphism.\end{cor}
\begin{prf}
We may without loss of generality assume that $S$ and $T$
are minimal.
Consider the Postnikov systems $S=S_\infty\rightarrow
\ldots\rightarrow
S_n\rightarrow \ldots\rightarrow S_1$ and $T=T_\infty\rightarrow\ldots
\rightarrow T_n\ldots\rightarrow T_1$ constructed above.
Clearly, all $S_i, T_i$ are minimal. In particular, $S_1$ and $T_1$ are 
weakly homotopy-equivalent minimal Kan complexes, thus are isomorphic (see
\cite{mayp}). 
This implies that $H_b^*\left(S_i\right)$ and $H_b^*\left(T_i\right)$
are isometrically isomorphic
for any $i$.

Finally, for any $n\in{\Bbb N}$, $S\rightarrow S_n$ induces an isomorphism
of $\pi_0,\ldots,\pi_n$ (\cite{mayp}).
By a simplicial version of the Hurewicz theorem (which {\em only in the case of Kan complexes} also holds true for bounded cohomology, here it would not suffice to assume the weak Kan property)
we get an isometric isomorphism $H_b^n\left(S_n\right)
\rightarrow H_b^n\left(S\right)$, the same way with $T_n$ and $T$. This 
proves $H_b^n\left(S\right)\simeq H_b^n\left(T\right)$ in (arbitrary)
degree n.
\end{prf}

\subsection{Moduli space of affine foliations}

In this section we study the question when two affine foliations of the standard simplex are topologically
conjugate.
For example, in dimension 2, there will be precisely the two non-conjugate affine foliations in the pictures
 below.
 
%\psset{unit=0.1\hsize}
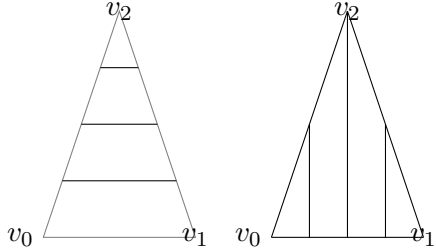
\begin{tikzpicture}
%\pspicture(0,0)(5,3)
\draw[gray](0,0)--(2,0)--(1,3)--(0,0);
\draw(0.5,1.5)--(1.5,1.5);
\draw(0.25,0.75)--(1.75,0.75);
\draw(0.75,2.25)--(1.25,2.25);
\node at(-0.3,0){$v_0$};
\node at(2,0){$v_1$};
\node at(1,3){$v_2$};

\draw(3,0)--(5,0);
\draw(5,0)--(4,3);
\draw(4,3)--(3,0);
\draw(3.5,0)--(3.5,1.5);
\draw(4,0)--(4,3);
\draw(4.5,0)--(4.5,1.5);
\node at(2.7,0){$v_0$};
\node at(5,0){$v_1$};
\node at(4,3){$v_2$};
\end{tikzpicture}

If a foliation is as in the right picture, we will say that the affine map has two isolated extrema (here $v
_0$ and $v_1$). If the foliation is as in the left picture,
we will say that it has
non-isolated extrema.\\
The result of this section which we are actually going to need is corollary 7, which
will be important for the proof of theorem 14.

\begin{lem}\label{affine} Let $f,g:\Delta^n\rightarrow{\Bbb R}$ be affine mappings such that
$$f\left(v_0\right)\ge f\left(v_1\right)\ge\ldots\ge f\left(v_n\right), g\left(v_0\right)\ge g\left(v_1
\right)
\ge \ldots\ge g\left(v_n\right)$$ and
such that for all $0\le i\le n-1$ we have $f\left(v_i\right)=f\left(v_{i+1}\right)$ if and only if $g\left(
v_i\right)=g\left(v_{i+1}\right)$. Then there is a piecewise linear homeomorphism $h:\Delta^n\rightarrow
\Delta^n$ such that for all $x,y\in\Delta^n$ we have: $$f\left(x\right)=f\left(y\right)\Longleftrightarrow g\left
(hx\right)=g\left(hy\right).$$\end{lem}
\begin{prf} Define $h\left(v_i\right)=v_i$ for $i=0,\ldots,n$.\\
We extend $h$ to the 1-skeleton. We will use that an affine mapping from a 1-simplex to ${\Bbb R}$ is either
 constant or injective.
Consider $\left[v_i,v_j\right]$ some 1-simplex with $i<j$. \\
If $f\left(v_i\right)=
f\left(v_j\right)$, then $f$ must be constant along $\left[v_i,v_j\right]$.
By assumption $g\left(v_i
\right)=g\left(v_j\right)$, thus also $g$ is constant along $\left[v_i,v_j\right]$. Then 
define $h$ to be the 
identity
on $\left[v_i,v_j\right]$. Thus $
f\left(x\right)
=f\left(y\right)\Longleftrightarrow g\left(hx\right)=g\left(hy\right)$ holds for
$x,y\in\left[v_i,v_j\right]$. \\
If $f\left(v_i\right)>f\left(v_j\right)$,
then $f$ is strictly increasing along $\left[v_i,v_j\right]$, hence there exist unique $
w_{i+1},\ldots,w_{j-1}\in\left[v_i,v_j\right]$ with $f\left(w_k\right)=f\left(v_k\right)$ for $i+1\le k\le j
-1$. The same way, since $g$ is 
strictly increasing (by assumption), we find unique
$u_{i+1},\ldots,u_{j-1}$ with $g\left(u_k\right)=g\left(v_k\right)$ for $i+1\le k\le j-1$. Define $h$ on $
\left[v_i,v_j\right]$ by $h\left(w_k\right)=u_k$
for $i+1\le k\le j-1$ and by piecewise linear extension, i.e., $$h\left(tw_k+\left(1-t\right)w_{k+1}\right)=
tu_k+\left(1-t\right)u_{k+1}$$ for all $t\in\left[0,1\right]$
and $i+1\le k\le j-1$. 

We are given $h$ on the 1-skeleton such that it is linear on $\left[w_k^{\left(i,j\right)},w_{k+1}^{\left(i,
j\right)}\right]\subset\left[v_i,v_j\right]$.
We want to extend $h$ to $\Delta^n$.
Consider $$W_i:=\left\{x\in\Delta^n:f\left(v_i\right)\ge
x\ge f\left(v_{i+1}\right)\right\}, U_i:=\left\{x\in\Delta^n:g\left(v_i\right)\ge x\ge g\left(v_{i+1}\right)
\right\}.$$
Since $f$ is affine, $W_i$
is a polytope with vertices $v_i, v_{i+1}, w_i^{\left(j,k\right)}$ for all
$j<i<k$ and $w_{i+1}^{\left(j,k\right)}$ for all $j<i+1<k$.
The same way, $U_i$ is a polytope
with vertices $v_i, v_{i+1}, u_i^{\left(j,k\right)}$ for all
$j<i<k$ and $u_{i+1}^{\left(j,k\right)}$ for all $j<i+1<k$.
We are given a bijection $h:w_i^{\left(j,k\right)}\rightarrow u_i^{\left(j,k\right)}$.
Fix some triangulation of $W_i$. Since $W_i$ and $U_i$ are combinatorially equivalent
polytopes, we may choose an 
equivalent triangulation for $U_i$. Moreover we may choose the triangulations of the $W_i$'s such that the triangulations
of $W_i$ and $W_{i+1}$ coincide on $W_i\cap W_{i+1}$ for all $i$. An affine mapping between simplices
is uniquely determined by their vertices.
Therefore we can define $h:W_i\rightarrow U_i$ as the unique piecewise linear mapping (with respect to the 
fixed triangulations) with $$h\left(v_i\right)=v_i, h\left(v_{i+1}\right)=v_{i+1},$$
$$ h\left(w_i^{\left(j,k\right)}\right)=u_i^{\left(j,k\right)}{\mbox\ for\ }j<i<k, h
\left(w_{i+1}^{\left(j,k
\right)}\right)=u_{i+1}^{\left(j,k\right)}{\mbox\ for\ }j<i+1<k.$$
We have $$h_i\mid_{V_i\cap V_{i+1}}=h_{i+1}\mid_{V_i\cap V_{i+1}}$$ because
the piecewise linear mapping is uniquely determined by its vertices. Thus
the $h_i$ fit together to a piecewise linear homeomorphism $h:\Delta^n\rightarrow\Delta^n$. We have:
$$f\left(x\right)=f\left(y\right)\Rightarrow x,y\in V_i\mbox{\ for\ some\ }i\Rightarrow g\left(hx\right)=g
\left(hy\right).$$
and vice versa.
$h$ is a continuous bijection between compact Hausdorff
spaces, thus a homeomorphism.\end{prf}\\
\\
{\bf Observation:} {\em Let $\left(M,{\mathcal{F}}\right)$ be a 
foliated manifold. Let $K$ be a simplicial complex 
built from copies $\Delta_1,\ldots,\Delta_r$ of the standard n-simplex $\Delta^n$ by identifications 
of faces. Assume that $\kappa:K\rightarrow M$ is a continuous mapping such that $\kappa_i=\kappa\mid_{\Delta_i}:\Delta^n=\Delta_i\rightarrow M$ is 
transverse to $\mathcal{F}$ for $i=1,\ldots,r$. Then there exists 
a simplicial homeomorphism $h:K\rightarrow K$ such that $\left(\kappa_ih\right)^*{\mathcal{F}}$ is
an affine foliation for $i=1,\ldots,r$.}\\
\begin{prf}
By transversality, we have 
that for each $\Delta_i$ such a simplicial homeomrphism $h_i:\Delta_i\rightarrow
\Delta_i$ exists, but a priori the $h_i$
need not agree on common faces. Therefore we construct $h$ on the $k$-skeleton of $K$, by induction
on $k$. Assume $h$ is defined on the $k-1$-skeleton, and $\tau$ is a $k$-dimensional subsimplex of $K$.
By transversality, there exists a simplicial homeomorphism $h_\tau:\tau\rightarrow\tau$ such that 
$h_{\tau}^*{\mathcal{F}}$ is an affine foliation. (If $k=1$, then $\tau$ need not be transverse, but of course ${\mathcal{F}}\mid_{\tau}$ is affine.)
W.l.o.g.\ $h_{\tau}$ maps each vertex (and hence each subsimplex of $\tau$) to itself.

On each face $\partial_j\tau$, the affine maps 
$f_{h_{\tau}}$ and $f_h$ (which define the affine foliations) satisfy, possibly after multiplication by 
$-1$, 
the assumptions of Lemma 6 (because $h_{\tau}$ and $h$ are simplicial), hence there exists $\overline{h}:
\partial\tau\rightarrow\partial\tau$ such that $f_{h_\tau}\circ \overline{h}=f_h$. 
($\overline{h}$ is constructed by induction on
the dimension of subsimplices of $\partial \tau$, to guarantee that it is well-defined on common faces of $\partial_{j_1}\tau$ and $\partial_{j_2}\tau$ for $j_1\not= j_2$.) In particular, $\left(h_{\tau}
\overline{h}\right)^*{\mathcal{F}}$ and $h^*{\mathcal{F}}$ are the same affine foliation on each face of $\tau$. 
%Thus exists a simplicial homeomorphism $\hat{h}$ of $\partial\tau$ with $\hat{h}h_{\tau}\overline{h}=h$ on $\partial\tau$. 
This implies that $\hat{h}:=h\mid_{\partial\tau}\left(h_{\tau}\mid_{\partial\tau}\overline{h}\right)^{-1}:\partial \tau\rightarrow
\partial\tau$ is a simplicial homeomorphism. 
Choose some arbitrary extensions of 
$\hat{h},\overline{h}$ as simplicial homeomorphisms of $\tau$. Then 
$h$ can be defined by $h:=\hat{h}h_{\tau}\overline{h}$ on $\tau$ and thus agrees with the definition of $h$ on $\partial\tau$.
\end{prf}\\

Referring to the title
of this section, \hyperref[affine]{Lemma \ref*{affine}} shows that the moduli space of affine foliations consists of a finite number of points.
More important for us will be the following corollary.

\begin{cor}\label{affi} Let $\Delta^n\subset {\Bbb R}^{n+1}$ be the standard simplex, $v_0,\ldots,v_n$ its vertices, and
 $v
$ some point in the interior of $\Delta^n$. Let $\tau_i$ be the straight simplex in ${\Bbb R}^{n+1}$
spanned by $v,v_0,\ldots,\hat{v}_i,\ldots,v_n$. Assume that affine mappings $f_i:\tau_i\rightarrow{\Bbb R}$
are defined such that the following conditions a),b),c),d) hold: \\
a) if $f_i\left(v_k\right)> f_i\left(v_l\right)$ for some $i,k,l$
%(and arbitrary vertices $u,w\in\left\{v,v_0,\ldots,\hat{v_i},\ldots,v_n$),
then $f_j\left(v_k
\right) >
f_j\left(v_l\right)$ for all $j$,\\
b) if $f_i\left(v_k\right)= f_i\left(v_l\right)$ for some $i,k,l$, then
$f_j\left(v_k
\right)=
f_j\left(v_l\right)$ for all $j$,\\
c) there do not exist $i,j,k,p,r,s$ such that \\
\hspace*{1in}$f_i\left(v_p\right)<f_i\left(v_r\right), f_j\left(v_r\right)<f_j\left(v_s
\right), f_k\left(v_s\right)< f_k\left(v_p\right)$\\
d) if we reindex $v_0,\ldots,v_n$ such that $f_i\left(v_k\right)\ge 
f_i\left(v_{k+1}\right)$ for all $i$ (this
is possible by a,b), then there is some $k$ with $f_i\left(v_k\right)\ge f_i\left(v\right)\ge f_i\left(
v_{k+1}\right)$ for all $i$.

Then there is an affine mapping $f:\Delta^n\rightarrow{\Bbb R}$ and a 
simplicial homeomorphism $h:\Delta^n\rightarrow 
\Delta^n$
such that $f\mid_{\tau_i}=f_ih$.\end{cor}
\begin{prf} We may define an ordering of $\left\{v_0,\ldots,v_n\right\}$ by $v_r\ge v_s\Leftrightarrow f_i
\left(v_r\right)\ge f_i\left(v_s\right)$ and $v_r>v_s\Leftrightarrow f_i\left(v_r\right)>f_i\left(v_s\right)$.
The assumption says that this is well-defined and transitive. Therefore there exists
some
affine mapping $f:\Delta^n\rightarrow{\Bbb R}$
with $f\left(v_r\right)>f\left(v_s\right)\Leftrightarrow v_r>v_s$ for all $1\le r,s\le n$. 
By d), we may choose
some $v^\prime\in
int\left(\Delta^n\right)$ such that $f\left(v_r\right)\ge f\left(v^\prime\right)\Leftrightarrow
 f_i
\left(v_r\right)\ge f_i\left(v^\prime\right)$ and $f\left(v_r\right) > f\left(v\right)\Leftrightarrow f_i\left(v_r
\right)
> f_i\left(v\right)$ for all $0\le r\le n$ and $i\not=r$.

For $0\le i\le n$ let $\tau_i^\prime\subset \Delta^n$ be the
straight simplex spanned by $v^\prime,v_0,\ldots,\hat{v_i},\ldots,v_n$. Let $g_i:\tau_i\rightarrow\tau_i^\prime
$ be the unique affine mapping which sends $v$ to $v^\prime$ and fixes $v_j$ for $j\not=i$.
 We are going to construct $h:\Delta^n\rightarrow
\Delta^n$.

$fg_i$ is affine. By \hyperref[affine]{Lemma \ref*{affine}}, we know that there exist homeomorphisms $h_i:\tau_i\rightarrow\tau_i$ such that
 $fg_i=f_ih_i$.
The homeomorphisms $h_i$ are constructed in the proof of lemma 6 in such a way
that they are uniquely defined as soon as triangulations of the different $U_k$'s and $W_k$'s are fixed. If
we choose the triangulations in such a way that they coincide on the common boundary faces $\tau_i\cap\tau_j
$, then
the so constructed homeomorphisms $h_i,h_j$ coincide on $\tau_i\cap\tau_j$:
$h_i\mid_{\tau_i\cap\tau_j}
=h_j\mid_{\tau_i\cap\tau_j}$. Therefore the $h_i$'s fit together to give a piecewiese linear homeomorphism $h
:\Delta^n\rightarrow\Delta^n$. $h$ fixes every face of $\Delta^n$, thus it is a simplicial homeomorphism.\\
Moreover, by construction, we have $g_i=g_j$ on $\tau_i\cap\tau_j$, therefore the $fg_i:\tau_i\rightarrow
{\Bbb R}$ fit together to a well-defined affine mapping.\end{prf}

\section{The norm on the fundamental group ...}
\subsection{... and its relation with the foliated Gromov norm}

\begin{defn}\label{fg} Let $\left(M,{\mathcal{F}}\right)$ be a foliated manifold,
$x_0\in M$ and $\gamma\in\pi_1\left(M,x_0\right)$. Then define \end{defn}
$$l_{\mathcal{F}}\left(\gamma\right)=inf\left\{\begin{array}{c}k: \mbox{\ exists representative $c$ of $\gamma$ and $k$ points $\left\{p_1,\ldots,p_k\right\}$
 on $c$ such that}\\
\ \ \ \ \ \ \ \mbox{each segment of $c-\left\{p_1,\ldots,p_k\right\}$ is strongly transverse to\ } {\mathcal{F}}\end{array}\right\}.$$
\\
This invariant may depend on the choice of base point $x_0\in M$ and may differ from the invariant $l\left(\alpha\right)$ which was defined in \cite{cale} as follows:
\begin{defn}\label{fg2} Let $\left(M,{\mathcal{F}}\right)$ be a foliated manifold and $\alpha:{\Bbb S}^1\rightarrow M$ a closed loop. Then define \end{defn}
$$l
\left(\alpha\right)=inf\left\{\begin{array}{c}k: \mbox{\ exists closed loop $\beta$ freely homotopic to $\alpha$ and $k$ 
points $\left\{p_1,\ldots,p_k\right\}$ on $\beta$}\\
\ \ \ \ \ \mbox{such that each segment of $\beta-\left\{p_1,\ldots,p_k\right\}
$ is strongly transverse to $\mathcal{F}$}\end{array}\right\}.$$
\\
There is an obvious inequality $l_{\mathcal{F}}\left(\left[\alpha\right]\right)\ge
l\left(\alpha\right)$, where $\left[\alpha\right]$ denotes the class of
$\alpha$ in $\pi_1\left(M,\alpha\left(0\right)\right)$.\\
It may happen that $l_{\mathcal{F}}\left(\left[\alpha\right]\right)$ and $l\left(\alpha\right)$ are not bounded. This is for example the case with the 
foliation discussed in \cite[Example 3.11]{cale}. On the other hand, classes of 
foliations with uniform bounds are exhibited in \cite{cale}. It is 
not clear what conditions on a foliation give uniform {\em nonzero} bounds on $l_{\mathcal{F}}$.\\
\\
We say that a codimension one foliation of an $m$-manifold $M$
satisfies the weak Kan property in degrees $\ge2$ if the union (over $n\in{\Bbb N}$) of
$$K_{n,st}^{\mathcal{F}}:=\left\{\sigma:\Delta^n\rightarrow M\mbox{\ 
singular simplex strongly 
transverse to\ }{\mathcal{F}}\right\}$$ satisfies the weak Kan property in degrees $2,3,\ldots,m=dim\left(M\right)$.\\
\\
Recall that $K_{st}^{\mathcal{F}}$ is a simplicial set. 
We observe that the inclusion $i:K_{*,st}^{\mathcal{F}}\rightarrow S_*\left(M\right)$ induces an {\bf isomorphism} $i_*:H_*^{simp}\left(K_{*,st}^{\mathcal{F}}\right)\rightarrow H_*\left(M\right)$. 
Indeed, let $\tau=\sum_{j=1}^r a_j\sigma_j$ be a singular chain with $\partial\tau=0$.
We may perform barycentrical subdivision (which preserves the homology class)
sufficiently often such that all simplices of the subdivision are
contained in a foliation chart and therefore can be homotoped to be strongly transverse.
If we perform
the homotopies succesively on simplices
of increasing dimensions, these homotopies are compatible (boundary
cancellations are preserved), because the homotopy may leave the $k-1$-skeleton 
pointwise fixed if the k-1-skeleton is already strongly transverse.
Hence the homotoped chain is a cycle in the homology class of $\tau$, showing that
$i_*$ is surjective. Similarly, if $\tau$ is strongly transverse and
$\tau=\partial\left(\sum_{j=1}^s b_j\kappa_j\right)$, we may perform barycentrical subdivision and homotopies on $\sum_{j=1}^s b_j\kappa_j$ to get a strongly transversal chain 
with boundary $\tau$, showing that $i_*$ is injective. \\
Since $i_*$ is an isomorphism, we have that, for any homology class $h\in H_*\left(M;{\Bbb R}\right)$, the equality $$\parallel h\parallel_{{\mathcal{F}}}^{st}=\parallel i_*^{-1}h\parallel_{simp}$$ holds, where $\parallel .\parallel_{simp}$ denotes the simplicial Gromov norm 
on $K_{*,st}^{\mathcal{F}}$. This will enable us to apply the results of Section 1.2.,
especially \hyperref[sap]{Proposition \ref*{sap}}. (By the way, a similar argument shows that 
$\parallel i_*^{-1}h\parallel_{sing}=\parallel h\parallel$, i.e. relations between $\parallel
h\parallel_{{\mathcal{F}}}^{st}$ and $\parallel h\parallel$ are actually relations between
$\parallel i_*^{-1}h\parallel_{simp}$ and $\parallel i_*^{-1}h\parallel_{sing}$.)\\
\\
The following theorem connects $l_{\mathcal{F}}$ with the foliated Gromov norm. The
technical assumption (of the weak Kan property being satisfied for degrees $\ge 2$)
is satisfied for fairly general foliations, as will be explained in Section 3. In particular, it is satisfied for any foliation without Reeb components on a 3-manifold.

\begin{thm}\label{ineq} Let $\left(M,{\mathcal{F}}\right)$ be a foliated manifold such that
$\widetilde{M}/\widetilde{\mathcal{F}}$ is a 1-connected 1-manifold
and the set of strongly transversal simplices satisfies the weak Kan property
in degrees $\ge 2$. Then $$\parallel h\parallel_{\mathcal{F}}\le
\left(1+\left(n+1\right)sup\left\{l_{\mathcal{F}}\left(\gamma\right):\gamma\in\pi_1
\left(M,x_0\right)\right\}\right)\parallel h\parallel$$
holds for any $x_0\in M$ and any $h\in H_n\left(M;{\Bbb R}\right)$.
In particular, if $M$ is closed and oriented,
$$\parallel M\parallel_{\mathcal{F}}\le
\left(1+\left(\dim\left(M\right)+1\right)sup
\left\{l_{\mathcal{F}}\left(\gamma\right)
:\gamma\in\pi_1
\left(M,x_0\right)\right\}\right)\parallel M\parallel.$$
\end{thm}

\begin{cor}\label{eq} Let $\left(M,{\mathcal{F}}\right)$ be a foliated, closed and oriented manifold such that
$\widetilde{M}/\widetilde{\mathcal{F}}$ is a 1-connected 1-manifold and
the set of strongly transversal simplices satisfies the weak Kan property
in degrees $\ge 2$. If $l_{\mathcal{F}}\left(\gamma\right)=0$
for some $x_0\in M$ and all $\gamma\in\pi_1\left(M,x_0\right)$, then
$$\parallel M\parallel_{\mathcal{F}}=\parallel M\parallel.$$\end{cor}

\begin{prf}
Let $\sum_{i=1}^r a_i\sigma_i$ be a cycle representing the homology class $h$.
We may homotope the $\sigma_i$ such that all vertices of all $\sigma_i$ are in the base point $x_0$. After this homotopy, all edges are closed loops $\gamma$, representing
classes $\left[\gamma\right]\in\pi_1\left(M,x_0\right)$, and we may further 
homotope, keeping $x_0$ fixed, such that the homotoped edges can be subdivided
into $l_{\mathcal{F}}\left(\left[\gamma\right]\right)$ transverse arcs (resp., if 
$l_{\mathcal{F}}\left(\left[\gamma\right]\right)=0$, such that the homotoped $\gamma$ is transverse to $\mathcal{F}$). It is straightforward to see that these homotopies can actually be extended to homotopies of $\sigma_1,\ldots,\sigma_r$ such 
that the homotopies do not affect cancellation of boundary faces. We continue to denote the homotoped singular
simplices by $\sigma_1,\ldots,\sigma_r$.

Let $L=sup\left\{l_{\mathcal{F}}\left(\gamma\right):\gamma\in\pi_1
\left(M,x_0\right)\right\}$.
We claim that we may subdivide each $\sigma_i$ into $1+\left(n+1\right)L$ simplices
$\tau_{ij}$ such that $\sigma_i$ is homologous to $\sum_j\tau_{ij}$ and
such that each $\tau_{ij}$ has transversal 1-skeleton.

Given a simplex $\sigma_i$, we can subdivide each of its edges $e$
by $l_{\mathcal{F}}\left(\left[e\right]\right)$ points into transversal arcs. Denote 
by $P$ this set 
of points. There is a (not unique) subtriangulation $T$
of $\sigma_i$ such that the vertices of simplices in $T$ are exactly the points in $P$ together with the vertices of $\sigma_i$. (We remark that these triangulations can be performed such that cancelling boundary faces of $\sigma_i$ and $\sigma_j$ 
are triangulated compatibly. This can just be achieved after prescribing an ordering on the total set of vertices.) The number of simplices in this triangulation is at most
$1+\left(n+1\right)L$. We claim that we can realize this subtriangulation of $\sigma_i$ such that its 1-skeleton is transverse.

To show this claim, 
we lift $\sigma_i$ to $\tilde{\sigma}_i$ in the universal covering $\widetilde{M}$. 
By assumption, the leaf
space $\widetilde{M}/\widetilde{\mathcal{F}}$ 
for the pull-back foliation $\mathcal
{F}$ is a simply connected 1-manifold. Therefore
the projection of $\tilde{\sigma}_i$ has (at least two) outermost vertices
in $\widetilde{M}/\widetilde{\mathcal{F}}$. Thus
it suffices to show the following: whenever an $n$-simplex $\tau$ is given with the property that for one of its outermost vertices 
$v\in\tau$ all edges emanating from $v$ are 
transverse, we can homotope $\tau$ to a simplex with transversal 1-skeleton.
(This implies the former claim, as we may apply the latter 
claim succesively to the simplices in 
the subtriangulation of $\sigma_i$, using that the 1-skeleton of $\sigma_i$ is
already transverse. At each step we have to take an outermost vertex of the remaining simplices 
in the subtriangulation $\sigma_i$, after the subsimplices whose 1-skeleton is already transverse 
have been removed.) So given an outermost vertex $v\in\tau$ and two other vertices $p_i,p_j$ of $\tau$, our task is to show that the arc connecting $p_i$ and $p_j$ can be homotoped to be transverse.

Consider the 
image of $\tau$ in $\widetilde{M}/\widetilde{\mathcal{F}}$. 
If $v$ is the outermost point, then, w.l.o.g., the projection of the
arc connecting $v$ to $p_j$ passes through $p_i$, hence the arc connecting $v$ and $p_j$ has to pass through the leaf $F_{p_i}$ containing $p_i$. Hence there is an arc connecting $p_i$ and $p_j$
which is composed by an arc in $F_{p_i}$ and by a transversal arc. A small perturbation makes this 
composed arc transverse to $\mathcal{F}$. (This argument shows also that, for any subtriangulation obtained after some removals, the vertex with the outermost projection to
$\widetilde{M}/\widetilde{\mathcal{F}}$ is an exterior vertex. Namely, if it was lying on an edge between two points $p_i$ and $p_j$, then the same argument would show that there is a transversal
arc connecting $p_i$ and $p_j$, i.e., we could reduce the number of points on this edge.)
This finishes the proof of the claim that we can realize the subtriangulation with transversal 1-skeleton.\\

Now we want to apply \hyperref[sap]{Proposition \ref*{sap}}.
We observe that $K^{\mathcal{F}}_{st}$ obviously is closed with respect to general degenerations, and that it also satisfies the simplicial approximation property. To see the latter statement, observe that each simplex can be barycentically subdivided sufficiently 
often such that all simplices of the iterated subdivision are contained in foliation charts. Inside these foliation charts, each simplex can be homotoped to a 
strongly transverse one. If we perform the homotopies succesively on simplices 
of increasing dimensions, these homotopies are compatible (boundary 
cancellations are preserved). This works because the homotopy may leave the $k-1$-skeleton pointwise fixed if it is already strongly transverse. 

Thus we have checked 
the assumptions of \hyperref[sap]{Proposition \ref*{sap}} and may apply \hyperref[sap]{Proposition \ref*{sap}} (with $K=K_{st}^{\mathcal{F}}$, which satisfies the
weak Kan property in degrees $\ge 2$). We have represented $i_*h$ by a chain $\sum_{i=1}^r a_i\sum_{j=1}^{1+\left(n+1\right)L}\tau_{ij}$ such that the 1-skeleta of all $\tau_{ij}$ belong to $K$. Thus we 
conclude $\parallel h\parallel_{\mathcal{F}}^{st}\le\left(1+\left(n+1\right)L\right)\parallel h\parallel$. This implies $\parallel h\parallel_{\mathcal{F}}\le\left(1+\left(n+1\right)L\right)\parallel h\parallel$. With $h=\left[M\right]$ we get the second claim of Theorem 1.
\end{prf}\\
\\
It should be mentioned that, for taut foliations of 3-manifolds, the implication\\
\hspace*{1in}$l\left(\alpha\right)=0\ \ \forall\alpha\Longrightarrow \parallel M\parallel_{\mathcal{F}}=\parallel M\parallel$\\
already follows from results in \cite{cale}.\\
\\
Remark: The proof of \hyperref[ineq]{Theorem \ref*{ineq}} would actually also work if we assumed that the set
of transversal simplices satisfies the weak Kan property in degrees $\ge2$. However, as
we will see in the beginning of section 3, there is {\em no} foliation  
${\mathcal{F}}$ such that $K_{\mathcal{F}}$ satisfies the weak Kan property in degree 2.
Thus, \hyperref[ineq]{Theorem \ref*{ineq}} would be (correct but) meaningless if we replaced `strongly transversal' by
`transversal' in the statement of the theorem.

\subsection {(Sub)Additivity of foliated Gromov norms}

Using \hyperref[sap2]{Lemma \ref*{sap2}}, one obtains the obvious generalization of \hyperref[ineq]{Theorem \ref*{ineq}} to foliated
manifolds with boundary. This is needed for the following corollary.
\begin{cor}\label{ineq2} Let $M_1,M_2$ be compact, orientable n-manifolds with foliations
${\mathcal{F}}_i, i=1,2$, transverse to
the connected boundaries and such that the leaf spaces of the universal coverings are 
1-connected 1-manifolds.
Assume that the boundaries are $\pi_1$-injective and
have amenable fundamental groups.\\
Let $f:\partial M_1\rightarrow \partial M_2$ a homeomorphism which is compatible
with ${\mathcal{F}}_1\mid_{\partial M_1}$ and ${\mathcal{F}}_2\mid_{\partial M_2}$,
 $M=M_1\cup_f M_2$ and $\mathcal{F}$ the glued foliation on $M$.
If all foliations satisfy the weak Kan property in degrees $n\ge 2$, then
$$\frac{1}{1+\left(n+1\right)L_1}\parallel M_1\parallel_{{\mathcal{F}}_1}+
\frac{1}{1+\left(n+1\right)L_2}\parallel M_2\parallel_{{\mathcal{F}}_2}\le$$
$$\ \ \ \ \ \parallel M\parallel_{\mathcal{F}}\le\left(1+\left(n+1\right)L\right)
\left(\parallel M_1\parallel_{{\mathcal{F}}_1}+
\parallel M_2\parallel_{{\mathcal{F}}_2}\right).$$
with $L=sup\left\{l_{\mathcal{F}}\left(\gamma\right):\gamma\in\pi_1
\left(M,x_0\right)\right\}$ and $L_i
=sup\left\{l_{{\mathcal{F}}_i}\left(\gamma\right):\gamma\in\pi_1
\left(M_i,x_0\right)\right\}$.

In particular, if all foliations satisfy the weak Kan property in degrees $n\ge 2$, and $l_{\mathcal{F}}\left(\gamma\right)=0$ for all $\gamma\in\pi_1\left(M,x_0\right)$
(for some $x_0\in M$), then 
$$\parallel M\parallel_{\mathcal{F}}=
\parallel M_1\parallel_{{\mathcal{F}}_1}+
\parallel M_2\parallel_{{\mathcal{F}}_2}.$$
\end{cor}
 \begin{prf}
The right hand inequality follows from \hyperref[ineq]{Theorem \ref*{ineq}} together with \\
$\parallel M\parallel=\parallel M_1,\partial M_1\parallel+\parallel M_2,\partial M_
2\parallel$ (\cite{grom},\cite{kw}) and the obvious inequality \\
$\parallel M_i,\partial M_i\parallel\le\parallel 
M_i,\partial M_i\parallel_{{\mathcal{F}}_i}$
for $i=1,2$.
Similarly one gets the left hand inequality.\end{prf}\\
\\
An analogous statement can be proved if the glueing only takes place along some connected components 
of $M_1$ resp.\ $M_2$, i.e., if $M$ happens to have nonempty boundary.\\
Moreover, if the foliated manifold $M_1$ has (at least) two (amenable, $\pi_1$-injective) boundary components
$\partial_1M_1,\partial_2M_1$ and one glues with a homeomorphism \\
$h:\partial_1M_1
\rightarrow\partial_2M_1$ to get a foliated manifold $\left(M=M_1/h,{\mathcal{F}}\right)$, one obtains \\
$\left(1+\left(n+1\right)L_1\right)^{-1}
\parallel M_1\parallel_{{\mathcal{F}}_1}\le\parallel M\parallel_{\mathcal{F}}
\le\left(1+\left(n+1\right)L\right)\parallel M_1\parallel_{{\mathcal{F}}_1}$\\
by the same arguments.\\
\\
\hyperref[ineq2]{Corollary \ref*{ineq2}} applies in particular to the JSJ-decomposition of 3-manifolds
carrying foliations without Reeb components.\\

\section{Weak Kan property for taut foliations}

Let $M$ be a compact manifold and $\mathcal{F}$ a codimension one foliation of $M$.

%We define a pseudonorm $\parallel h\parallel^{strong}$ on homology by admitting only strongly transversal simplices. Clearly
%$\parallel h\parallel^{strong}\ge\parallel h\parallel$.

Denote $$ K_n^{\mathcal{F}}:=\left\{\sigma:\Delta^n\rightarrow M\mbox{\ singular simplex transverse to\ }
{\mathcal{F}}\right\}$$
and
$$K_{n,st}^{\mathcal{F}}:=\left\{\sigma:\Delta^n\rightarrow M\mbox{\ singular simplex strongly transverse to\ }{\mathcal{F}}\right\}.$$
Recall that the union
$K_{st}^{\mathcal{F}}:=\cup_{n\ge 0}K_{n,st}^{\mathcal{F}}$ 
is stable with respect to boundary maps and degeneracy maps,
meanwhile, in general, the faces of simplices in $K_2^{\mathcal{F}}$ 
do not necessarily belong to $K_1^{\mathcal{F}}$. (This is due to the fact that each foliated
{\em degenerate} 2-simplex belongs to $K_2^{\mathcal{F}}$, but not necessarily to $K_{2,st}^{\mathcal{F}}$.)\\
We mention that, 
for an arbitrary foliation $\mathcal{F}$, the simplicial set $K_{st}^{\mathcal{F}}$ does not satisfy the weak Kan property in degree 1, and $K^{\mathcal{F}}$ does not satisfy the weak Kan property neither in degree 1 nor in degree 2.
To prove the first statement, just consider two transversal 1-simplices whose 
composition is not transverse. To prove the latter statement, consider (for the foliation of ${\Bbb R}^n$ by horizontal hyperplanes with the induced total ordering of 
the leaf space) 4 points $v_0,v_1,v_2,v_3$ such that, with respect to the total ordering, $v_0>v_1>v_2=v_3$ (for the leaves containing the respective vertices) and consider the straight simplices
$\sigma_2$ resp.\ $\sigma_3$ with vertices $v_0,v_1,v_3$ resp.\ $v_0,v_1,v_2$. Let $\sigma_0$ be the degenerate 2-simplex with edges $\left(v_1,v_2\right)$ and $\left(
v_1,v_3\right)$. Although $\sigma_0$ is transverse (but not strongly transverse)
to $\mathcal{F}$, it is clear that there is no transverse 2-simplex $\sigma_1$
with $\partial_0\sigma_1=\partial_0\sigma_0,\partial_1\sigma_1=\partial_0\sigma_2,
\partial_2\sigma_1=\partial_0\sigma_3$. Hence we get a contradiction to 
the weak Kan property. 
Clearly, the same argument works for any foliated manifold $\left(M,{\mathcal{F}}\right)$
because it can be realised inside a foliation chart.

\begin{lem}\label{weak}Let $M^3$ be a compact 3-manifold, $\mathcal{F}$ a codimension one foliation on $M$. Then \\
a) the simplicial sets $K_{st}^{\mathcal{F}}$ and $K^{\mathcal{F}}$
satisfy the weak Kan property in degree 3, \\
b) the simplicial set $K^{\mathcal{F}}$ 
does not satisfy the weak Kan property in degree 2, \\
c) the simplicial set $K_{st}^{\mathcal{F}}$ satisfies
the weak Kan property in degree 2 if and only if $\mathcal{F}$ has no Reeb component.\\
\end{lem}

We will deduce \hyperref[weak]{Lemma \ref*{weak}} from the following, more general, \hyperref[thm2]{Theorem \ref*{thm2}}.\\
We say that a foliation $\mathcal{F}$ of a manifold is {\bf fibration-covered}
if there is some covering $\widehat{M}\rightarrow M$ with 
pull-back foliation $\widehat{\mathcal{F}}$ on 
$\widehat{M}$, such that the projection $\widehat{\pi}:\widehat{M}\rightarrow\widehat{M}/\widehat{\mathcal{F}}$ from $\widehat{M}$ to its leaf space is a locally trivial fibration. 
(An equivalent condition is that just
$\tilde{\pi}:\widetilde{M}\rightarrow\widetilde{M}/\widetilde{{\mathcal{F}}}$ is a locally trivial fibration, where $\widetilde{M}$ denotes the universal covering and $\widetilde{{\mathcal{F}}}$ the pull-back foliation.)
This property is more common than one might expect in view of the following observation:\\
\\
{\bf Observation.} 
If $M$ is an $m$-dimensional 
manifold and $\mathcal{F}$ a codimension 
one foliation satisfying the following conditions:\\
(i) for every leaf $F$ is $\pi_1F\rightarrow\pi_1M$ injective,\\
(ii) for every leaf $F$, the universal covering 
$\widetilde{F}$ is homeomorphic to ${\Bbb R}^{m-1}$,\\
then  $\mathcal{F}$ is fibration-covered and $\widetilde{M}/\widetilde{\mathcal{F}}$ is a 1-connected 1-manifold.\\
\\
\begin{prf} Let $\widetilde{M}$ be the universal covering with the pull-back foliation
$\widetilde{\mathcal{F}}$.
The assumptions imply that all leaves of $\widetilde{\mathcal{F}}$ are homeomorphic to ${\Bbb R}^{m-1}$. These foliations have been investigated in \cite{palm} and it has been shown, in particular,
that the projection $\tilde{\pi}:\widetilde{
M}\rightarrow \widetilde{M}/\widetilde{{\mathcal{F}}}$
is a locally trivial fibration 
(this is the corollary to the trivialization lemma on \cite{palm},
p.117). Moreover, in chapter 1 of \cite{palm} it is shown that $\widetilde{M}/\widetilde{\mathcal{F}}$ is a 
1-connected 1-manifold.
\end{prf}

\begin{thm}\label{thm2}
If $M$ is an $m$-dimensional
manifold and $\mathcal{F}$ a  fibration-covered codimension
one foliation, such that $\widetilde{M}/\widetilde{\mathcal{F}}$ is a 1-connected 1-manifold, then 

a) $K_{st}^{\mathcal{F}}$ satisfies 
the weak Kan property in degrees $2,\ldots,m$, and

b) $K^{\mathcal{F}}$ satisfies the weak Kan property in degrees $3,\ldots,m$.\end{thm}
\begin{cor}\label{weak2}
If $M$ is an $m$-dimensional manifold and $\mathcal{F}$ a codimension one 
foliation satisfying the following conditions:\\
(i) for every leaf $F$, $\pi_1F\rightarrow\pi_1M$ is injective,\\
(ii) for every leaf $F$, the universal covering $\tilde{F}$ 
is homeomorphic to 
${\Bbb R}^{m-1}$,\\
then

a) $K_{st}^{\mathcal{F}}$ satisfies
the weak Kan property in degrees $2,\ldots,m$, and

b) $K^{\mathcal{F}}$ satisfies the weak Kan property in degrees $3,\ldots,m$.\end{cor}

\begin{prf}

Consider n-simplices $\tau_0,\ldots,\tau_n:\Delta^n\rightarrow M$ with $\partial_i\tau_j=\partial_{j-1}\tau_i$ for all $i<j$. 
Let $K$ be the simplicial complex built from $n+1$ copies $\Delta_0,\ldots,\Delta_n$
of the 
standard simplex $\Delta^n$ by identifying $\partial_i\Delta_j$ with $\partial_{j-1}\Delta_i$ for all $i<j$.
It comes with a map $\tau:K\rightarrow M$, which on $\Delta_i$ agrees with $\tau_i$.
$K$ is isomorphic to the union $\partial_0\Delta^{n+1}\cup\ldots\cup\partial_n\Delta^{n+1}$ of all but one faces of the standard n+1-simplex $\Delta^{n+1}$. We will denote $v:=\partial_0^n\Delta_i\in K$ the common vertex of all $\Delta_i$ and $v_i\in K$ the unique vertex of $K$ which is not contained in $\Delta_i$. For simplicity, we 
will write $\tau_i$ for the mapping $\tau_i:\Delta^n\rightarrow{\Bbb R}$ as well
as for the subsimplex $\Delta_i\subset K$.
For each $\tau_i$ we have that there exists a simplicial homeomorphism $h_i:\Delta^n\rightarrow\Delta^n$
such that $\left(\tau_i h_i\right)^*{\mathcal{F}}$ 
is an affine foliation. By the observation after lemma 6, this implies that
there exists a simplicial homeomorphism $h:K\rightarrow K$ such that $\left(\tau_i h\right)^*{\mathcal{F}}$ is an affine foliation of each $\tau_i$. Therefore we will henceforth assume that $\tau_i^*{\mathcal{F}}$ is an affine foliation for $i=0,\ldots,n$.\\

Of course, whenever $n$-simplices $\tau_0,\ldots,\tau_n$ satisfy
$\partial_i\tau_j=\partial_{j-1}\tau_i$ for all $i<j$, we have the canonical degenerate $n+1$-simplex $\sigma$ with $\partial_i\sigma=\tau_i$ for all $i\le n$ which is defined by precomposing $\tau_0\cup\ldots\cup\tau_n$ with the canonical retraction from the standard $n+1$-simplex $\Delta^{n+1}$
to the union of $\partial_0\Delta^{n+1},\ldots,\partial_n\Delta^{n+1}$. We will denote this simplex by $\sigma_{deg}\left(
\tau_0,\ldots,\tau_n\right)$. (The cases of interest will be those where $\partial_{n+1}\sigma_{deg}\left(\tau_0,\ldots,\tau_n\right)$
is not transverse to $\mathcal{F}$ even though $\tau_0,\ldots,\tau_n$ are.)

We assume, for $i=1,\ldots,n$, that $\tau_i^*{\mathcal{F}}$ is a foliation of $\Delta_i$ homeomorphic to a foliation by level sets of an affine mapping
$f_i:\Delta_i\rightarrow{\Bbb R}, 0\le i\le n$.

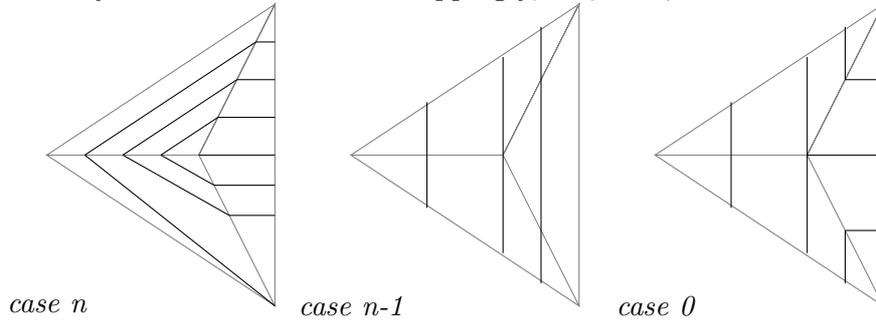
\begin{tikzpicture}
\draw[gray](0,2)--(3,4)--(2,2)--(0,2)--(3,0)--(2,2)--(3,4)--(3,0);
\draw(1.5,2)--(2.25,2.5);
\draw(1,2)--(2.5,3);
\draw(0.5,2)--(2.75,3.5);
\draw(2.2,1.6)--(1.5,2);
\draw(2.4,1.2)--(1,2);
\draw(3,0)--(0.5,2);
\draw(2.25,2.5)--(3,2.5);
\draw(2.5,3)--(3,3);
\draw(2.75,3.5)--(3,3.5);
\draw(2,2)--(3,2);
\draw(2.2,1.6)--(3,1.6);
\draw(2.4,1.2)--(3,1.2);
\node at(0,0){\em case n};

\draw[gray](4,2)--(7,4)--(6,2)--(4,2)--(7,0)--(6,2)--(7,4)--(7,0);
\draw(5,1.3)--(5,2.7);
\draw(6,0.7)--(6,3.3);
\draw(6.5,0.3)--(6.5,3.7);
\node at(4,0){\em case n-1};

\draw[gray](8,2)--(11,4)--(10,2)--(8,2)--(11,0)--(10,2)--(11,4)--(11,0);
\draw(9,1.3)--(9,2.7);
\draw(10,0.7)--(10,3.3);
\draw(10.5,0.3)--(10.5,1);
\draw(10.5,1)--(11,1);
\draw(10,2)--(11,2);
\draw(10.5,3)--(10.5,3.7);
\draw(10.5,3)--(11,3);
\node at(8,0){\em case 0};

\end{tikzpicture}

We will distinguish two possibilities. Possiblity $A$ is that there exists some simplex $\tau_i$ such that $v$ is not an extremum of $f_i$. (Some examples are pictured above, where case $k$ means that $v$ is an extremum of $k$ simplices.)
In this case we will show, without needing any assumption on the foliation $\mathcal{F}$, that $\partial_{n+1}\sigma_{deg}\left(\tau_0,\ldots,\tau_n\right)$ is strongly transverse to $\mathcal{F}$ if $\tau_0,\ldots,\tau_n$ are. (If $n\ge 3$, we will also get transversality of $\partial_{n+1}\sigma_{deg}\left(\tau_0,\ldots,\tau_n\right)$ under 
the weaker assumption that $\tau_0,\ldots,\tau_n$ are transverse, not necessarily strongly transverse.)
Possibility $B$ is that $v$ is an extremum of all $f_i$. In this case we will use the assumptions on $\mathcal{F}$ to construct a simplex $\tau_{n+1}$ (strongly) transverse to $\mathcal{F}$.\\

{\bf A: $v$ is not an extremum of $f_0$.}\\

We claim that in this case there exists some vertex which is an extremum of all of $f_0,\ldots,f_n$ (except for the simplex which does not contain this vertex, of course).\\
First note that there are $n+1$ simplices with (at least) $2n+2$ extrema, but only $n+2$ vertices. If $v$ is
not an extremum of all $f_0,\ldots,f_n$, then some $v_k\not=v$ must be the extremum of at least two different simplices. W.l.o.g. we may assume that $v_0$ is an extremum of $f_1,f_2$.

We \underline{claim} that this implies that \\
- either \underline{\bf $v_0$ is an extremum of $f_1,\ldots,f_n$}\\
- or $v_0$ belongs to the same leaf as $v,v_3,\ldots, v_n$ and \underline{$v_1,v_2$ are extrema of all simplices} (where they occur). \\
For $n=2$, there is nothing to prove, so we will restrict to the assumption $n\ge 3$.\\
To prove the claim by contradiction, assume that $v_0$ is not an extremum of $f_3$. Then we have: $$f_3\left(v_k\right)>f_3\left(v_0\right)>f_3\left(v_l\right)$$ for some $v_k,v_l$.
Since we are assuming $n\ge 3$, at least one of the following three cases holds: \\
- there exists $v_m\not=v_k$ with $f_3\left(v_k\right)\ge f_3\left(v_m\right)>f_3\left(v_0\right)>f_3\left(v_l\right)$,\\
- or there exists $v_m\not=v_l$ with $f_3\left(v_k\right)>f_3\left(v_0\right)>f_3\left(v_m\right)\ge f_3\left(v_l\right)$,\\
- or $f_3\left(v_k\right)>f_3\left(v_0\right)=\ldots=f_3\left(v_j\right)>f_3\left(v_l\right)$.\\
Note that the first two cases are equivalent after replacing $f_3$ with $-f_3$, so it suffices to consider one of them.\\
Consider case 1, i.e., $f_3\left(v_k\right)\ge f_3\left(v_m\right)>f_3\left(v_0\right)>f_3\left(v_l\right)$. 
If $l\not =1$ and $m\not=1$, then we have 
$v_m>v_0>v_l$ on $\tau_1\cap\tau_3$, i.e., $v_0$ is not an extremum of $f_1\mid_{\tau_1\cap\tau_3}$, contradicting the fact that it is an extremum of $f_1$. 
If $l\not=1$ and $k\not=1$, then we have $v_k>v_0>v_l$ on $\tau_1\cap\tau_3$, i.e., $v_0$ is not an extremum of $f\mid_{\tau_1\cap\tau_3}$, contradicting the fact that it is an 
extremum of $f_1$. Since $k$ and $m$ can not be both equal to $1$, we have derived a contradiction if $l\not=1$.\\
If $l=1$ and $m\not=2$, then $v_m>v_0>v_1$ in $\tau_2\cap\tau_3$, i.e., 
$v_0$ is not an extremum of $f_2\mid_{\tau_2\cap\tau_3}$,
contradicting the fact that it is an 
extremum of $f_2$. If $l=1$ and $k\not=2$, then $v_k>v_0>v_1$ in $\tau_2\cap\tau_3$, i.e., $v_0$ is not an extremum of $f_2\mid_{\tau_2\cap\tau_3}$,
contradicting the fact that it is an 
extremum of $f_2$. Thus we have derived a contradiction also if $l=1$. This finishes case 1 (and the equivalent case 2.)\\
In case 3, we conclude that all vertices except $v_k,v_l$ and possibly $v_3$
belong to the same leaf.
If $k,l\not=1$, we get a contradiction because $v_0$ would not be an extremum of $f_3\mid_{\tau_1\cap\tau_3}$. If $k,l\not=2$,
we get a contradiction because $v_0$ would not be an extremum of $f_3\mid_{\tau_2\cap\tau_3}$.
There remain the cases $k=1,l=2$ resp.\ $l=1,k=2$ which need some more care. If $k=1,l=2$, then
$$f_3\left(v_1\right)>f_3\left(v_0\right)=f_3\left(v\right)>f_3\left(v_2\right).$$
We distinguish the two possibilities that $v_3$ belongs to the same leaf as $v,v_0,v_4,\ldots,v_n$ or not. \\
Consider first the case that $v_3$ does not belong to the same leaf as $v$ and $v_0$.
Consider the subtetrahedron $T_0$ of the simplex $\tau_0$ which is spanned by the 
four vertices $v_1,v_2,v_3,v$. Since $v$ is not an extremum of $f_3\mid_{T_0\cap\tau_3}$ 
but is an extremum of $f_1\mid_{T_0\cap\tau_1}$, and 
since $v_1,v_2$ are extrema of $f_3\mid_{T_0\cap\tau_3}$, 
we necessarily have either $$f_0\left(v_1\right)>f_0\left(v\right)>f_0\left(v_3\right)>f_0\left(v_2\right)$$
or $$f_0
\left(v_1\right)>f_0\left(v\right)>f_0\left(v_2\right)>f_0\left(v_3\right)$$
(possibly after replacing $f_0$ with $-f_0$).\\
In both cases we have $v_1>v>v_3$ on $T_0\cap\tau_2$, hence necessarily
$$f_2\left(v_1\right)>f_2\left(v_0\right)=f_2\left(v\right)>f_2\left(v_3\right)$$
(possibly after replacing $f_2$ by $-f_2$) because $v_0$ and $v$ belong to the same leaf.
But this contradicts the assumption that $v_0$ is an extremum of $f_2$.\\
Thus we are left with the case that all vertices except $v_1$ and $v_2$ belong to the 
same leaf. Recall that $v_1$ and $v_2$ are isolated extrema of $f_3$. Let $m\not=1,2,3.$ Looking at $\tau_m\cap\tau_3$, we note
that $v_1$ and $v_2$ are extrema of $f_3\mid_{\tau_m\cap\tau_3}$. 
Since $f_m\left(v_3\right)=f_m\left(v_i\right)$ for $i\not=1,2$, this implies 
$$f_m\left(v_1\right)> f_m\left(v\right)=f_m\left(v_0\right)=f_m\left(v_4\right)=\ldots f_m\left(v_n\right)>f_m\left(v_2\right).$$
Finally, for $m=1,2$ we trivially have (possibly after replacing $f_1$ by $-f_1$ or $f_2$ by $-f_2$): $$f_2\left(v_1\right)>f_2\left(v\right)=f_2\left(v_0\right)=f_2\left(v_3\right)=\ldots=f_2\left(v_n\right),$$
$$f_1\left(v\right)=f_1\left(v_0\right)=f_1\left(v_3\right)=\ldots=
f_1\left(v_n\right)> f_1\left(v_2\right).$$
Hence, $v_1$ and $v_2$ are extrema of all $f_i$.\\
If $k=2,l=1$, the same argument works.\\

So, we now can
assume that some vertex, say $v_0$, is an
extremum of all affine mappings (except $f_0$, where it does not occur).
Replacing some $f_i$ by $-f_i$ if necessary, we have that $v_0$ is a maximum of $f_1,\ldots,f_n$. We \underline{claim} that this implies that 
the \\
\underline{\bf assumptions of \hyperref[affi]{Corollary \ref*{affi}} are satisfied} for $\partial_{n+1}\sigma_{deg}\left(\tau_0,\ldots,\tau_n\right)$.

First, since $v_0$ is a maximum of $f_1,\ldots,f_n$ we have that $$f_j\left(u\right)>f_j\left(w\right)\Longleftrightarrow f_i\left(u\right)>f_i\left(w\right),
f_j\left(u\right)=
f_j\left(w\right)\Longleftrightarrow f_i\left(u\right)=f_i\left(w\right)$$
whenever $i,j\not=0$ and $u,w\in\left\{v_1,\ldots,v_n,v\right\}-\left\{v_i,v_j\right\}$. Indeed, if not, we would get a contradiction by looking at the induced foliation of the triangle $\Delta\left(v_0uw\right)$ in $\tau_i\cap\tau_j$.
(Hence we have checked the assumptions a,b,c of corollary 7 except for $f_0$.)

In particular, there exists some $u\in\left\{v_1,\ldots,v_n,v\right\}$ such that $f_i\left(u\right)\ge f_i\left(w\right)$ for all $1\le i\le n$
and $w\in\left\{v_1,\ldots,v_n,v\right\}-\left\{u,v_i\right\}$. This vertex $u$
is then a maximum of $f_i\mid_{\tau_0\cap\tau_i}$ for all $1\le i\le n$ except possibly for $\tau_0\cap\tau_j$ if $u=v_j$.\\

For the rest of the proof we have to distinguish the cases $n\ge 3$ and $n=2$.\\
Consider $n\ge 3$. In this case, $u$ being a maximum of $f_i\mid_{\tau_0\cap\tau_i}$ for all $1\le i\le n$ implies
that $u$ is a maximum of $f_0\mid_{\tau_0\cap\tau_i}$ for all $1\le i\le n$ and therefore
that $u$ is a maximum of $f_0$. Then we may argue as before: $f_0\left(w\right)>f_0\left(w^\prime\right)$ is equivalent to
$f_i\left(w\right)>f_i\left(w^\prime\right)$ for any $i$ because, if not, we would get a contradiction in $\Delta\left(uww^\prime\right)\subset \tau_0\cap\tau_i$. Hence we get conditions a,b,d for \hyperref[affi]{Corollary \ref*{affi}}.

 It remains to check condition c of \hyperref[affi]{Corollary \ref*{affi}}.\\
Trivially, $f_i\left(v_p\right)<f_i\left(v_r\right)$ and $f_j\left(v_r\right)<f_j\left(v_s\right)$ together imply $f_k\left(v_j\right)> f_k\left(v_p\right)$ except possibly if $j=p$ or $i=s$.
If $j=p$, then we have to check that $f_i\left(v_j\right)<f_i\left(v_r\right),
f_j\left(v_r\right)<f_j\left(v_s\right)$ and $f_k\left(v_s\right)<f_k\left(v_j\right)$ can not 
happen simultaneously. However, these three inequalities would lead, looking at the triangle $\Delta\left(v_jv_rv_s\right)$ with the induced foliation, to the contradiction $v_j<v_r<v_s<v_j$. To be precise, choose some $q\not=i,j,r$ (this is possible if $n\ge 3$), then we have $f_q\left(v_j\right)<f_q\left(v_r\right)<f_q\left(v_s\right)<f_q\left(v_j\right)$, getting a contradiction.\\
If $i=s$, an analogous argument works. This finishes the proof for $n\ge 3$.

We are left with the case $n=2$.\\
The conditions a,b,c of \hyperref[affi]{Corollary \ref*{affi}} are empty for $n=2$, we want to check condition d.\\
Case 1: $u=v$. (This will be the step in the proof which uses that we are working with strongly transverse 
simplices rather then just transverse ones.) We claim that $v$ is an extremum 
of $f_0$. If not, then either the leaf through 
$v$ would intersect $K$ in a trivalent graph (if 
$\tau_0$ is nondegenerate) as in the picture of {\em case n} above, which is of course impossible for the leaf being a manifold, or $\tau_0$ would be degenerate. 
We claim that in the latter case, $\tau_0$ could not be strongly transverse. Indeed, looking at a small 
foliation chart around $v$, we observe that the leaf space of this 
chart (which is just an open intervall) is decomposed by the leaf $L_v$ through $v$ into two components, and that the image of $\tau_0$ in the leaf space
is completely contained in the closure of one of these components. This means that the map from $\tau_0$ to the 
leaf space is not locally surjective 
at $\tau_0^{-1}\left(L_v\right)$, contradicting the strong transversality of $\tau_0$. Thus, $v$ is an extremum of $f_0$. Together with $f_1\left(v_0\right)\ge f_1\left(v\right)\ge f_1\left(v_2\right)$ and $f_2\left(v_0\right)\ge f_2\left(v\right)\ge f_2\left(v_1\right)$ this implies that condition d holds true.\\
Case 2: $u=v_1$. This means $f_2\left(v_0\right)\ge f_2\left(v_1\right)\ge f_2\left(v\right)$.
We know that $v_0$ is a maximum of $f_1$. Thus one possibilty for $f_1$ is $f_1\left(v_0\right)\ge f_1\left(v\right)\ge f_1\left(v_2\right)$. 
Then, if $v$ were an extremum of $f_0$, we would be in {\em case n} pictured at the beginning of the proof (with 
degenerate $\tau_1$). Clearly, $\tau_1$ would not be strongly transverse, giving a contradiction. Thus, $v$ must 
not be an extremum of $f_0$, 
from which condition d of corollary 7 easily follows.
The other possibility for $f_1$ is $f_1\left(v_0\right)\ge f_1\left(v_2\right)
\ge f_1\left(v\right)$. But then it is immediate 
that $\tau_0$ is not strongly transverse, hence this possibility 
can not happen.\\
The case $u=v_2$ works the same way. \\

We have checked the assumptions of \hyperref[affi]{Corollary \ref*{affi}}.
Thus we may apply \hyperref[affi]{Corollary \ref*{affi}} to get an affine mapping
$f:\partial_{n+1}\sigma_{deg}\left(\tau_0,\ldots,\tau_n\right)\rightarrow{\Bbb R}$. This shows transversality. 
(We note that we have used the strong transversality of $\tau_0,\ldots,\tau_n$ only in one step, namely in the
special discussion to handle the case $n=2$. Thus, for $n\ge 3$, we get part a) of \hyperref[thm2]{Theorem \ref*{thm2}}.)

We continue with the proof for $K^{\mathcal{F}}_{st}$, i.e.\ we wish to prove that
$\partial_{n+1}\sigma_{deg}\left(\tau_0,\ldots,\tau_n\right)$ is even strongly transverse. We assume that $\tau_0,\ldots,\tau_n$ are strongly transverse. Thus, the only points where strong transversality of 
$\partial_{n+1}\sigma_{deg}\left(\tau_0,\ldots,\tau_n\right)$
may fail are points in some intersection $\tau_i\cap\tau_j$, where necessarily
this intersection must belong to the same leaf. However, looking at the foliation chart
in a neighborhood of $v$ we would get that the leaf space is decomposed by the leaf through $v$ into two components,and that any other simplex as $\tau_i$ and $\tau_j$ can therefore not be strongly transverse at $v$,
giving a contradiction. (One may look at the picture for {\em case n} for understanding.) Thus,
$\partial_{n+1}\sigma_{deg}\left(\tau_0,\ldots,\tau_n\right)$
is automatically strongly transverse.\\

\begin{tikzpicture}

\draw[gray](0,2)--(3,4)--(2,2)--(0,2)--(3,0)--(2,2)--(3,4)--(3,0);
\draw(1.5,2)--(2.25,2.5);
\draw(1,2)--(2.5,3);
\draw(0.5,2)--(2.75,3.5);
\draw(2.25,2.5)--(2.2,1.6);
\draw(2.5,3)--(2.4,1.2);
\draw(2.75,3.5)--(3,0);
\draw(2.2,1.6)--(1.5,2);
\draw(2.4,1.2)--(1,2);
\draw(3,0)--(0.5,2);
\draw(0.2,2)--(2.9,3.8);
\draw(2.9,3.8)--(3,2);
\draw(0.2,2)--(1.5,1);
\node at(0,0){\em case n+1a};

\draw[gray](4,2)--(7,4)--(6,2)--(4,2)--(7,0)--(6,2)--(7,4)--(7,0);
\draw(5.5,2)--(6.25,2.5);
\draw(6.25,2.5)--(6.2,1.6);
\draw(6.2,1.6)--(5,2);
\draw(5,2)--(6.5,3);
\draw(6.5,3)--(6.4,1.2);
\draw(6.4,1.2)--(4.5,2);
\draw(4.5,2)--(6.75,3.5);
\node at(4,0){\em case n+1b};

\end{tikzpicture}

{\bf B: $v$ is an extremum of $f_0,\ldots,f_n$.}\\
\\
We assume that $f_i:\tau_i\rightarrow{\Bbb R}, i=0,\ldots,n$ are affine
mappings such that the
common vertex $v:=\partial_o^n\tau_i$ of the n+1 $n$-simplices is an 
isolated extremum of $f_0,\ldots,f_n$, as pictured above. (We will see that a picture as in case n+1b actually can not happen.)
Recall that we have a simplicial complex $K$ and a map $\tau:K\rightarrow M$.
Since $K$ is contractible, we may lift $\tau$ to $\tilde{\tau}:K\rightarrow \widetilde{M}$. 
$K$ inherits the foliation $\tilde{\tau}^{-1}\left(\widetilde{\mathcal{F}}\right)$. \\
Each leaf $\widetilde{F}$ of
$\widetilde{\mathcal{F}}$ separates $\widetilde{M}$, since $\widetilde{M}/\widetilde{\mathcal{F}}$ is a 
1-connected 1-manifold. 
%(If $\widetilde{F}$ would not 
%separate, one could construct a loop $l$ intersecting $\widetilde{F}$ only 
%once, contradicting $\pi_1\widetilde{M}=0$. Indeed, in the 1-point compactification $\widetilde{M}^+$, the homological intersection
%number of $l$ and $\widetilde{F}^+$ would be 1, but $l$ is 0-homologous.) 
% $Clearly, $\pi_1\left(\widetilde{M}\right)\rightarrow\pi_1\left(\widetilde{M}/\widetilde{{\mathcal{F}}}\right)$ is surjective (any loop 
%can be lifted). In particular, $\pi_1\left(\widetilde{M}/\widetilde{\mathcal{F}}\right)=0$, i.e., 
%the leaf space of $\widetilde{\mathcal{F}}$ is simply connected. This implies, since $\dim\left(
%\widetilde{M}/\widetilde{\mathcal{F}}\right)=1$, that every point in the leaf space of $\widetilde{\mathcal{F}}$ separates, hence every leaf $\widetilde{F}$ of $\widetilde{\mathcal{F}}$ separates.\\
This means in particular that no connected component of $\tilde{\tau}^{-1}\left(\widetilde{F}\right)$
can intersect an edge of some $K$ twice, that is, a 
situation as in the picture `case n+1b' above does not happen. Namely, if 
$\tilde{\tau}^{-1}\left(\widetilde{F}\right)$
intersected an edge twice,
then the image of this edge under $\tilde{\tau}$
would connect points in the two path-components of $\widetilde{M}-\widetilde{F}$, giving the contradiction that $\widetilde{M}-\widetilde{F}$
were path-connected.

As a consequence we have that, for each leaf $F$ of $\mathcal{F}$,
the leaves $\tilde{\tau}_i^{-1}\left(\widetilde{F}\right)=\tilde{\tau}_i^{-1}\left(F\right)\subset K$,
which are (up to homeomorphism) n+1 
$(n-1)$-dimensional simplices, fit together (as in the picture of case n+1a) in a standard way, that is: \\
- close to the common vertex $v$
they fit together to $(n-1)$-dimensional continuous spheres,\\
- there may be leaves which do not intersect all simplices, but they fit together to proper subsets of $(n-1)$-dimensional continuous spheres.\\
In other words, the foliation of $K$ is (up to homeomorphism)
precisely the foliation that one would obtain if $K$ were to be embedded into the ${\Bbb R}^n$, foliated by 
horizontal planes, such that the embedding is linear on each $\tau_i$.

Let $\tilde{\pi}:\widetilde{M}\rightarrow \widetilde{M}/\widetilde{\mathcal{F}}$ be 
the projection. Strong transversality of $\tau_i:\Delta^n\rightarrow M$ implies that the image
of $\tilde{\pi}\tilde{\tau}_i:\Delta^n\rightarrow \widetilde{M}/\widetilde{\mathcal{F}}$ is an (unbranched) interval. Since the induced foliations on $\tau_i,\tau_j\subset K$ agree, the images of $\tilde{\pi}\tilde{\tau}_i$ for distinct $i$ are (subsets of) the same interval $I$. 

We are assuming that $\tilde{\pi}:\widetilde{M}\rightarrow
\widetilde{M}/\widetilde{\mathcal{F}}$ is a locally trivial fibration. Since $I$ is paracompact and contractible, the restriction of $\tilde{\pi}$ to $\tilde{\pi}^{-1}\left(I\right)$
is a trivial fibration, i.e.,
$$\tilde{\pi}^{-1}\left(I\right)\simeq I\times 
\widetilde{F}$$ is a product. This product structure makes the construction of the $n$-simplex with transversal $n+1$-th boundary face kind of obvious. Namely,
using this product structure, the foliation of 
$K=\tau_0\cup\ldots\cup\tau_n$ corresponds to a continuous family of continuously embedded $n-1$-spheres in $\widetilde{F}$, parametrized by a connected subinterval of $I$, at least 
until the first vertex $\not= v$ is reached. After the first vertex the 
family continues as a family of proper subsets of spheres. The picture is, as mentioned before, the same that one would get by embedding $K$ piecewise linearly into the standard foliation. \\
Observe that all spheres represent trivial elements in $\pi_{n-1}\left(\widetilde{F}\right)$, since the starting point of the continuous family has been a sphere mapped to a point. (The parametrization by a connected subinterval of $I$ is giving the homotopy.)
This means that this continuous family
of continuously embedded spheres in $\widetilde{F}$ can 
be extended to a continuous family of continuously embedded $n$-balls 
in the respective leaves, parametrized over the same subinterval of $I$).

If $v_0,\ldots,v_n$ belong to the same leaf $F_0$, then 
this continuous family of balls allows us to define a continuous mapping $\sigma$ from
the standard $(n+1)$-simplex to $\widetilde{M}$ such that $\partial_{n+1}\sigma$ is contained in the leaf $F_0$, thus is strongly transverse.

If $v_0,\ldots,v_n$ do not belong to the same leaf, then we get a 
continuous family of balls until the first vertex $\not=v$, say $v_0$, 
is reached. Let $F_0$ be the leaf containing $v_0$ and assume $v_n\not\in F_0$. 
We may use this family of balls to homotope $\tau_0,\ldots,\tau_n$ 
by homotoping $v$ (and leaving $\partial_n\tau_0,\ldots,\partial_n\tau_n$
fixed) until $v\in F_0$, preserving (strong) transversality.

After this homotopy we are in a situation that $v_0$ is a
common extremum for the induced foliations
on the homotoped simplices $\hat{\tau_0},\ldots,\hat{\tau_n}$. We have shown in the course of the proof of case A that
this implies that the assumptions of \hyperref[affi]{Corollary \ref*{affi}} are satisfied for 
$\tau_{n+1}:=\partial_{n+1}\sigma_{deg}\left(\hat{\tau_0},\ldots,\hat{\tau_n}\right)$. Hence,
by \hyperref[affi]{Corollary \ref*{affi}},  $\tau_{n+1}
$ is transverse and, as we have seen in the proof of case A, then automatically strongly
transverse since $\hat{\tau}_0,\ldots,\hat{\tau}_n$ are strongly transverse. Since $\partial_i\tau_{n+1}=\partial_n\hat{\tau}_i
=\partial_n\tau_i$, this finishes the proof.

\end{prf}\\
\\

We are now ready to prove \hyperref[weak]{Lemma \ref*{weak}}.\\
\begin{prf}
First note that b) holds for any foliation $\mathcal{F}$, as
has been observed at the beginning of Section 3.\\
We discuss a) and c).
In the course of the proof, we will use some deep results from foliation theory. Assume that $\left(
M,{\mathcal{F}}\right)$ is
a foliated 3-manifold. One distinguishes the case that there exists a leaf homeomorphic to the 2
-sphere or not. If {\bf some leaf
is homeomorphic to ${\Bbb S}^2$}, the Reeb stability theorem implies that
the foliation is a locally trivial fibration with fiber ${\Bbb S}^2$.
In this case, arguments analogous to the proof of theorem 14 imply that the weak Kan property in
degree 2 is fulfilled for $K^{\mathcal{F}}_{st}$. (The proof of case A did not use any assumptions on $\mathcal{F}$, and in the discussion of case B, it suffices to use that the fibers are simply
connected to fill 1-spheres by 2-balls, thus getting the weak Kan property in degree 2.)
To prove the weak Kan property in degree 3, for
${\Bbb S}^2$-bundles over ${\Bbb S}^1$, one observes that, even though the fibers are not 2-connected, the argument in possiblity B still works. Namely, working in the universal covering ${\Bbb S}^2\times{\Bbb R}$, one observes that,
 after an identification of all fibers with a fixed fiber, we have a continuous family of 2-spheres terminating in a constant 2-sphere (coming from the lift
of the fiber through $v$), in particular all 2-spheres arising from the intersection are 0-homotopic in their respective fibers, thus can be fibre-wise
filled with 3-balls. This allows again to define a transversal 4-simplex with transversal boundary faces.\\
\\
 If there is {\bf no leaf homeomorphic to ${\Bbb S}^2$} and the foliation
has {\bf no Reeb component} (i.e., no compressible torus leaf), it follows
from Palmeira's work that the induced foliation of the universal cover $\left(\widetilde{M},
\widetilde{\mathcal{F}}\right)$
is a foliation of $\widetilde{M}={\Bbb R}^3$ by leaves homeomorphic to ${\Bbb R}^2$ such that 
each leaf separates
${\Bbb R}^3$ into two connected components. (Namely, if there is no Reeb component,
then all leaves $F$ are $\pi_1$-injective, by Novikov's theorem. Hence 
$\widetilde{F}\simeq{\Bbb R}^2$, i.e. $\widetilde{\mathcal{F}}$ is
a foliation by planes. Then apply \cite[Corollary 3]{palm}.) 
It is shown in \cite{palm} that such a foliation satisfies both assumptions of \hyperref[thm2]{Theorem \ref*{thm2}}.
Thus we may apply \hyperref[thm2]{Theorem \ref*{thm2}} to get the weak Kan property in 
degrees 2 (for $K^{\mathcal{F}}_{st}$) and 3 
(for $K^{\mathcal{F}}$ and $K^{\mathcal{F}}_{st}$).\\
\\
It remains to discuss the case that {\bf $\mathcal{F}$ has Reeb components} (which will actually not be used in the paper)
and to show that $\mathcal
{F}$ nevertheless satisfies the
weak Kan property in degree 3, but not in degree 2. \\

We start with discussing degree 2.
We want to show that existence of a Reeb component destroys the weak Kan property 
(for the set of strongly transversal simplices) in degree 2.
Let $\mathcal{F}$ be the Reeb foliation of ${\Bbb D}^2\times {\Bbb S}^1$. Let
$\tau_0,\tau_1,\tau_2$ be 3 triangles,
whose images have a common vertex in 
$\left(0,\alpha\right)\in {\Bbb D}^2\times {\Bbb S}^1$,
and which satisfy 
$\partial_0\tau_0=\partial_0\tau_1, \partial_0\tau_2=\partial_1\tau_0
, \partial_1\tau_1=\partial_1\tau_2$. Assume that the images of
$\partial_2\tau_0, \partial_2\tau_1,
 \partial_2\tau_2$ lie on the boundary torus 
$\partial \left({\Bbb D}^2\times 
{\Bbb S}^1\right)$. Let $K$ be the simplicial complex built from 3 copies $\Delta_0,\Delta_1,\Delta_2$ of the standard simplex $\Delta^2$ by identifying $\partial_0\Delta_0=\partial_0\Delta_1, 
\partial_0\Delta_2=\partial_1\Delta_0,
\partial_1\Delta_1=\partial_1\Delta_2$, and let $\tau:K\rightarrow 
{\Bbb D}^2\times {\Bbb S}^1$ be the continuous mapping with $\tau\mid_{\Delta_i}=\tau_i$ for $i=0,1,2$.
Then $\tau\mid_{\partial K}:\partial K\rightarrow
{\Bbb D}^2\times {\Bbb S}^1$
represents a generator of $$ker\left(\pi_1\left(\partial\left(
{\Bbb D}^2\times{\Bbb S}^1\right)\right)\rightarrow\pi_1\left({\Bbb D}^2\times {\Bbb S}^1
\right)\right)\simeq{\Bbb Z}.$$
We can choose the 2-simplices $\tau_0,\tau_1,\tau_2$
to be strongly transverse to the Reeb foliation.
For any 3-simplex $T$ with $\partial_iT=\tau_i$ for $i=0,1,2$ we have that $\tau_3:=
\partial_3T$ is a 2-simplex
with $\partial_0\tau_3=\partial_2\tau_0,\partial_1\tau_3=\partial_2\tau_1,\partial_2
\tau_3=\partial_2\tau_2$. In particular, ${\mathcal{F}}
\mid_{\tau_3}$
contains $\partial\tau_3$ as a leaf, because
$\partial\left({\Bbb D}^2\times{\Bbb S}^1\right)$ 
is a leaf of $\mathcal{F}$. Thus ${\mathcal{F}}
\mid_{\tau_3}$ is a foliation of a topological disk such that the boundary is a leaf.
 Any such foliation of a topological disk must have a singularity (because the Euler
 characteristic of the disk does not vanish), in particular can not be topologically conjugate to
 an affine foliation of the
2-simplex. Thus, $\partial_3T$ can not be transverse, except if its image
were contained in the torus leaf $\partial\left({\Bbb D}^2\times{\Bbb S}^1\right)$.
%But the latter is impossible because, by construction, 
However, $\tau:K\rightarrow 
{\Bbb D}^2\times {\Bbb S}^1$
%the union $im\left(\tau_0
%\right)\cup im\left(\tau_1\right)\cup im\left(\tau_2\right)$
must represent a nontrivial element in $$\pi_2\left({\Bbb D}^2\times{\Bbb S}^1,\partial
\left({\Bbb D}^2\times{\Bbb S}^1\right)\right)\simeq{\Bbb Z},$$ 
because its boundary (in the long exact homotopy sequence) is $\tau\mid_{\partial K}:\partial K\rightarrow 
{\Bbb D}^2\times {\Bbb S}^1$, which was nontrivial.
Hence, if $\tau_3$ 
were contained in $\partial\left({\Bbb D}^2\times{\Bbb S}^1\right)$,
there could be no 3-simplex $T$ with $\partial_iT=\tau_i$ for $i=0,1,2,3$. Therefore 
the image of $\partial_3T$ can not be contained in
${\Bbb D}^2\times {\Bbb S}^1$, thus $\partial_3T$ can not be transverse. This shows that 
a fliation with a Reeb component does not satisfy the weak Kan property in degree 2.\\

Finally we show that any foliation $\mathcal{F}$
of a 3-manifold $M$ satisfies the weak Kan property in degree 3. 
We have already seen this for foliations without Reeb components.\\
Let $R\not=\emptyset$ be the union of the Reeb components. We remark that the leaf space of $\widetilde{\mathcal{F}}\mid_{\widetilde{int\left(R\right)}}$ is ${\Bbb R}$ and the leaf space of 
$\widetilde{\mathcal{F}}\mid_{\widetilde{\overline{M-R}}}$ is a 1-connected 1-manifold because ${\mathcal{F}}\mid_{\overline{M-R}}$ is a Reebless foliation.

Let transversal 3-simplices $\tau_1,\tau_2,\tau_3,
\tau_4$ be given. 
Let $K$ be the simplicial complex built from 4 copies $\Delta_0,\Delta_1,\Delta_2,\Delta_3$ 
of the standard simplex $\Delta^3$ by identifying
$\partial_{j-1}\Delta_i=\partial_i\Delta_j$ for all $i<j$.
Let $\tau:K\rightarrow M$
be the continuous mapping with $\tau\mid_{\Delta_i}=\tau_i$ for $i=0,1,2,3$. Let $\partial K:=\partial_0\Delta_1\cup\partial_0\Delta_2
\cup\partial_0\Delta_3\cup\partial_0\Delta_4$.
Let $x_0$ be the common vertex of $\Delta_1,\ldots,\Delta_4$ and let $x_i$ be the vertex not contained in 
$\Delta_i$, for $i=1,\ldots,4$. Let $v_i=\tau\left(x_i\right)$ for $i=0,\ldots,4$.
We distinguish the following 4 cases:\\
a) $v_0\in int\left(R\right), v_1,\ldots,v_4\not\in int\left(R\right)$,\\
b) $v_0\not\in int\left(R\right), v_1,\ldots,v_4\in int\left(R\right)$,\\
c) $v_0\not\in int\left(R\right), v_i
\not\in int\left(R\right)$ for some $i\in\left\{1,\ldots,4\right\}$,\\
d) $v_0\in int\left(R\right), v_i\in int\left(R\right)$ 
for some $i\in\left\{1,\ldots,4\right\}$.

{\em Case a)}: We observe that in this case $\tau\left(\partial K\right)$
can not intersect $int\left(R\right)$. Indeed, if it did, then some edge $\left[v_i,v_j\right]:=
\tau\left(\left[x_i,x_j\right]\right)$ 
would have to intersect $int\left(R\right)$ (because, if 
a leaf intersects a transversal simplex, then it must intersect 
its 1-skeleton). Hence there would be some subintervall $\left[w_i,w_j\right]\subset\left[v_i,v_j\right]$
contained in $R$, such that $w_i,w_j\in\partial R$. We note that the leaf space of the 
Reeb foliation of any connected component
of $int\left(R\right)$ is homeomorphic to ${\Bbb R}$. Hence the image of $\left[w_i,w_j\right]$ in the leaf space has some extremal point. But at this extremal point, $\left[w_i,w_j\right]$ would not be a submersion, hence not stronly transverse.

Since $\tau\left(\partial K\right)$ does not intersect $int\left(R\right)$ and since
$$\pi_3\left(R,\partial R\right)=0$$ we may homotope $\tau\left(K\right)$ off $int\left(R\right)$,
leaving $\tau\mid_{\partial K}$ fixed. This 
homotopy can be made simplicial. Let $\hat{\tau}:K\rightarrow M$
%\hat{\tau_0}\cup\hat{\tau_1}\cup\hat{\tau_2}\cup\hat{\tau_3}$
be the result of the homotopy, and let $\hat{\tau}_i:
\Delta_i\rightarrow M$ be defined by $\hat{\tau}_i=\hat{\tau}\mid_{\Delta_i}$.
Since the image of $\hat{\tau}$ does not intersect
$int\left(R\right)$ (and we have just proved that the weak Kan property 
holds in degree 3 for foliations without Reeb components), we do
get a 4-simplex $\hat{L}\in S_*\left(M- int\left(R\right)\right)$ such that
$\partial_i\hat{L}=\hat{\tau_i}$ for $i=0,1,2,3$ and
such that $\partial_4\hat{L}$ is transverse to $\mathcal{F}$. Since the Kan property is, of course,
true for $S_*\left(M\right)$, we can use $\hat{L}$ and the simplicial homotopy
to
produce, by successive application of the Kan property, a (not necessarily transversal)
4-simplex $L$ with
$\partial_4 L=\partial_4\hat{L}$ and $\partial_i L=\tau_i$ for $i=0,1,2,3$. 
Since $\partial_4\hat{L}$ is transverse, this finishes the proof in case a).

{\em Case b)}: This is similar to case a). We observe that the image of $\partial L$ can not intersect 
$N:=M - int\left(R\right)$. Indeed, if it did, we would again find some subinterval
$\left[w_i,w_j\right]\subset \left[v_i,v_j\right]$ contained in $N$ with $v_i,v_j\in\partial N$. Let ${\mathcal{G}}:={\mathcal{F}}\mid_N$ and $\widetilde{\mathcal{G}}$ its pull-back to
the universal covering $\widetilde{N}$. The leaf space of the Reebless foliation
$\widetilde{\mathcal{G}}$ 
is a 1-connected 1-manifold. Hence the image of $\left[w_i,w_j\right]
$ in the leaf space would have
some extremal point, giving a contradiction. 

It is well-known that the existence of a Reebless foliation implies that $N$ is irreducible and has infinite 
fundamental group. Hence $\pi_3\left(N\right)=0$. Since $\partial N$ consists of tori, 
this implies $\pi_3\left(N,\partial N\right)=0$. Thus we can homotope $L$ into
$int\left(R\right)$ and then apply the same argument as in case a). (Note 
that the restriction of ${\mathcal{F}}$ to $int\left(R\right)$ is a trivial foliation by planes, for which the weak Kan property in degrees $\ge 2$ holds.)

{\em Case c)}: If $v_0$ is not a common extremum of all ${\mathcal{F}}\mid_{\tau_i}$, we are in the situation of case A in the proof of \hyperref[thm2]{Theorem \ref*{thm2}}. We know then from the proof of that case that, without any assumptions on ${\mathcal{F}}$, there exists a
4-simplex $L$ with $\partial_iL=\tau_i, i=1,2,3,4$, such that $\partial_0L$ is transverse to $\mathcal{F}$.

So we are left to consider the case that $v_0$ is a common extremum of ${\mathcal{F}}
\mid_{\tau_i}, i=1,2,3,4$. We claim that we can homotope 
$\tau$ to some $\hat{\tau}$,
%=\hat{\tau}_1\cup\hat{\tau}_2\cup\hat{\tau}_3\cup\hat{\tau}_4$, 
leaving $\tau\mid_{\partial K}$ fixed, such that $v_0$ is 
not a common extremum of ${\mathcal{F}}\mid_{\hat{\tau}_i}$ for all $i=1,2,3,4$. After 
having accomplished this homotopy, we can apply the above argument to $\hat{\tau}$ instead of $\tau$.

The homotopy can be accomplished by an argument similar to that in case B in the proof of \hyperref[thm2]{Theorem \ref*{thm2}}. Let $N=M -int\left(R\right)$. ${\mathcal{G}}={\mathcal{F}}\mid_N$
is a Reebless foliation, hence we know from the proof of case B that $\tilde{\pi}^{-1}\left(I\right)=I\times \widetilde{G}$ for a connected interval $I\subset \widetilde{N}/\widetilde{\mathcal{G}}$ and a leaf $\widetilde{G}$
of $\widetilde{\mathcal{G}}$. Let $v_i$ be the vertex contained in $N$,
which exists by assumption of case c). We remark that $v_i$
must belong to the same connected component of $N$ as $v_0$. (Otherwise the edge 
$\left[v_0,v_i\right]$ would enter and leave 
some Reeb component, hence could not be transverse.)
Now we can copy the argument which we used in the proof of case B in \hyperref[thm2]{Theorem \ref*{thm2}}.
Namely, we can again use the product structure, to get a continuous family of 3-balls until
the first vertex $\not= v_0$ is reached. (We may assume w.l.o.g.\ that $v_i$ is this
first vertex.) Let $F_0$ be the leaf containing $v_i$.
We may use the family of balls to homotope $\tau_1,\ldots,\tau_4$
by homotoping the image of $x_0$ (and leaving $\partial_0\tau_1,\ldots,\partial_0\tau_4$
fixed) until the image of $x_0$ belongs to $F_0$, preserving strong transversality.
After this homotopy we are in a situation that the image of
$x_0$ is an extremum for the induced foliations
on the homotoped simplices $\hat{\tau_1},\ldots,\hat{\tau_4}$. We have shown in the course of
the proof of case A of \hyperref[thm2]{Theorem \ref*{thm2}}
that
this implies that the assumptions of \hyperref[affi]{Corollary \ref*{affi}} are satisfied for $\tau_0:=\partial_0
\sigma_{deg}\left(\hat{\tau_1},\ldots,\hat{\tau_4}\right)$. Hence $\tau_0
$ is strongly transverse. Since $\partial_i\tau_0=\partial_0\hat{\tau}_i
=\partial_0\tau_i$, this finishes the proof.

{\em Case d)}: 
If $v_0$ and $v_i$ belong to the same connected component of $R$, then the argument is literally the
same as in case c), with the roles of $N$ and $R$ interchanged.
If $v_0\in R_1$ and $v_i\in R_2$ for distinct connected components $R_1,R_2$ of $R$,
then we can apply the same argument as in case a) to homotope $K$ off $int\left(R_1\right)$.
If there is some $v_j\not\in int\left(R\right)$, then we are, after this homotopy, in the situation of case c) and can finish the proof. So it remains
to discuss the case that $v_1,v_2,v_3,v_4\in int\left(R\right)$.

We assume that $v_0$ is a common extremum of $\tau_1,\ldots,\tau_4$. (Else we are done by the proof of case A in \hyperref[thm2]{Theorem \ref*{thm2}}.)
We can apply the argument of case c) to homotope 
the image of $x_0$, 
leaving $\partial_0\tau_1,\ldots,\partial_0\tau_4$ fixed and preserving strong transversality, 
until the image of $x_0$ is contained in the boundary of a Reeb component $R_j$ with 
$v_j\in R_j$ for some $j\in\left\{1,2,3,4\right\}$. Since $v_0$ is a common extremum,
it follows that
all four points $v_1,v_2,v_3,v_4$ must belong to this Reeb component $R_j$. By the same argument as in case b), this implies that no edge $\left[v_i,v_l\right]$ can leave $R_j$,
for $i,l\in\left\{1,\ldots,4\right\}$. 
In particular, $\partial R_j$ does not intersect any edge $\left[v_i,v_l\right]$ and hence, by transversality, no 2-simplex $\left(v_iv_lv_k\right)$ with $i,l,k\in\left\{1,2,3,4\right\}$.
That is, these four 2-simplices are contained in $int\left(R_j\right)$. Since these four 2-simplices are strongly transverse and the foliation of $int\left(R_j\right)$ is a trivial
product foliation, we 
can find a strongly transversal 3-simplex, whose boundary faces are these four 2-simplices. 
This finishes the proof.

\end{prf}\\
\\
{\bf Generalizations.} In \cite{cale}, the notion of foliated Gromov norm has 
been generalized to laminations and, more generally, to group actions on order trees. We briefly describe the analogous generalization of the transversal length on the
fundamental group.

Let $M$ be a manifold and $\Gamma=\pi_1M$ its fundamental group, which acts 
on the universal covering $\widetilde{M}$. Let $\Gamma$ act on an order tree $T$ (see \cite{goe} for the definition of order tree) and assume that there
is a $\Gamma$-equivariant map $\phi:\widetilde{M}\rightarrow T$. (Such an 
$\phi$ exists naturally for any lamination $\mathcal{F}$ of $M$.) According to 
\cite[Definition 4.2.1]{cale} a singular i-simplex $\sigma:\Delta^i\rightarrow M$
is transverse if, for any lift $\widetilde{\sigma}:\Delta^i\rightarrow\widetilde{M}$, the image of $\phi\widetilde{\sigma}$ is a totally ordered segment of $T$. We say that a singular $i$-simplex is strongly transverse if it is transverse and the induced mapping $\phi\tilde{\sigma}:\Delta^i\rightarrow T$ is a submersion, for any lift $\tilde{\sigma}$ of $\sigma$ to $\widetilde{M}$.\\
Then one can again define the notions of transversal length on the
fundamental group and foliated Gromov norm.\\
It is easy to see that the proof of \hyperref[ineq]{Theorem \ref*{ineq}} also works in this setting.
However, to get useful information from the generalized
\hyperref[ineq]{Theorem \ref*{ineq}}, it would be necessary to know under what conditions the simplices strongly transverse to a given lamination (resp.\ group action on an order tree) satisfy the weak Kan property in degrees $\ge2$.

\end{document}